\documentclass[12pt]{article}
\usepackage{amssymb,amsmath,amsthm,epsfig}
\usepackage{latexsym, enumerate}
\usepackage{eepic}
\usepackage{epic}
\usepackage{graphicx}
\usepackage{color}
\usepackage{ifpdf}
\usepackage{subfigure}

\usepackage{dsfont}
\usepackage{multirow}
\usepackage{makecell}
\usepackage{algorithm}

\theoremstyle{plain}

\theoremstyle{definition}

\theoremstyle{remark}
\newtheorem{rem}{Remark}[section]

\begin{document}
\title{ \large\bf Model's sparse representation based on reduced  mixed GMsFE basis methods}

\author{
Lijian Jiang\thanks{Institute  of Mathematics, Hunan University, Changsha 410082, China. Email: ljjiang@hnu.edu.cn. Corresponding author}
\and
Qiuqi Li\thanks{College of Mathematics and Econometrics, Hunan University, Changsha 410082, China. Email:qiuqili@hnu.edu.cn.}
}

\date{}
\maketitle

\begin{center}{\bf ABSTRACT}
\end{center} %\smallskip

In this paper, we propose a model's sparse representation based on reduced mixed generalized multiscale finite element (GMsFE) basis methods for elliptic PDEs with random inputs.
A typical application for the elliptic PDEs is the flow in heterogeneous random porous media.
Mixed generalized multiscale finite element method (GMsFEM) is one of the accurate and efficient approaches to solve the flow problem in a coarse grid and obtain the velocity with local mass conservation.
When the inputs of the PDEs are parameterized by the random variables, the  GMsFE basis functions usually depend on the random parameters.  This leads to  a large number  degree of freedoms for the mixed GMsFEM and
substantially impacts on the computation efficiency.  In order to overcome the difficulty,  we develop reduced mixed GMsFE basis methods such that the multiscale basis functions are independent of the random parameters
and span a low-dimensional space.  To this end, a greedy algorithm is used to find a set of  optimal samples from a training set scattered in the parameter space.  Reduced mixed GMsFE basis functions  are constructed based on the optimal samples using two optimal sampling strategies: basis-oriented cross-validation and proper orthogonal decomposition.  Although the dimension of the space spanned by the reduced  mixed GMsFE basis functions is much smaller than the dimension of the original full order model, the online computation still depends on the number of coarse degree of freedoms.  To significantly improve the online computation,
we integrate the reduced mixed GMsFE basis methods with sparse tensor approximation and obtain a sparse representation for the model's outputs.  The sparse representation is very efficient
for evaluating the model's outputs for many instances of parameters.  To    illustrate the efficacy of the proposed methods, we present a few numerical examples for elliptic PDEs with multsicale and random inputs.
In particular,  a two-phase flow model in random porous media is simulated by the proposed   sparse representation method.

{\bf keywords:} Mixed generalized multiscale finite element method,  Basis-oriented cross-validation method,   Least-squares method of snapshots, Sparse tensor approximation

\section{Introduction}

Many complex real-world problems of physical and engineering interests have multiple scales and  uncertainties, for example, the subsurface flow models in highly heterogeneous porous media. Due to the measurement noise and the lack of knowledge about the physical properties, the model inputs usually contain some uncertainties. The uncertainties are often
parameterized by random variables. Thus, parameterized partial differential equations (PPDEs) are used to characterized the models with uncertainties.
One of the great challenges in these models is to  efficiently and accurately solve the PPDEs with high-contrast multiple scales  to predict the model's outputs and estimate the model's parameters. Both the multiple  scales and uncertainties have great challenges  for  numerical simulation. To reduce computation complexity and improve simulation efficiency, model reduction methods are desirable to the complex multiscale models.
In recent years, many numerical  methods have been proposed  to solve the complex  PPDEs using  multiscale methods and model reduction methods \cite{clc11, pde08, jpl12, ncn08,rhp08, aej08, arbogast2007multiscale, efendiev2013generalized, jiang2010stochastic}.

In many practical applications such as multiphase flow in porous media,   local mass conservation  is necessary for a numerical method to solve the problems.
The mixed finite element methods retain local conservation of mass and have been found to be
useful for solving these problems.
To illustrate the main idea of our approach, we consider the following mixed formulation of an elliptic  PDE with random inputs,
\begin{eqnarray}
\label{mix-eq}
\begin{cases}
\begin{split}
 k(x,\mu)^{-1}v(\mu)+\nabla p(\mu)&=0 ~~\text{in} ~~D, \\
\nabla\cdot v(\mu)&=f(x) ~~\text{in}~~D,
\end{split}
\end{cases}
\end{eqnarray}
where $D \in R^2$ is a bounded spatial domain, $f(x)$ is a source term, and $k(x,\mu)$ is a random permeability field, which may be oscillating with respect to the random parameter $\mu$ and  highly heterogeneous  with respect to the spatial variable $x$.  To take account of the heterogeneities, multiscale methods  \cite{abgr05, arbogast2007multiscale, ee03, hughes98,  pde08, efendiev2009multiscale}
can be used to solve the problem.  In this work, our interest is to use mixed multiscale finite element methods \cite{chen2003mixed, cel15}.
 The main idea of mixed multiscale finite element methods is to compute   multiscale basis functions for each interface supported in the coarse blocks  sharing the common interface.
 In order to accurately capture complex multiscale features, a few multiscale basis functions may be required for each coarse block. This is the idea of  Generalized Multiscale Finite Element Method (GMsFEM) \cite{cel15, efendiev2013generalized}.
 We note that the multiscale basis functions can be computed overhead and used repeatedly for the model with different source terms, boundary conditions and the coefficients with similar multiscale structures \cite{efendiev2009multiscale}.

Because the coefficient $k(x,\mu)$ in equation (\ref{mix-eq}) varies  with the realizations of random parameter $\mu$, the GMsFE basis functions usually depend on the parameters.
This significantly effects on the computation efficiency. In order to get the multiscale basis functions independent of the parameter, we build multiscale basis functions based on a set of samples in the parameter space. This will result in a high-dimensional GMsFE  space for approximation and bring great challenge for simulation. To alleviate the high dimensionality of the GMsFEM space, we
 identify a set of optimal reduced basis functions from the high dimensional GMsFE space and develop reduced mixed GMsFE basis methods. Thus
  a reduced order multiscale model can be obtained by projecting the original full order model onto the reduced multiscale finite element space.

In this work, we focus on  mixed GMsFEM and present reduced  mixed GMsFE basis methods to solve parameterized elliptic PDEs. Reduced basis  method is one of effective model order reduction methods and has been used to solve PPDEs in a low-dimensional manifold \cite{ buffa2012priori, chen2010certified, ncn08, rhp08}. The main idea of the reduced basis  method is to construct a small set of basis functions based on a set of snapshots, which are the solutions of the PPDEs corresponding to a set of optimal parameter samples selected by some sampling strategies. An offline-online computational decomposition is achieved  to improve efficiency. In offline stage, snapshots are computed and reduced basis functions are generated. In online stage,  the output is computed by the resultant reduced model for many instances of parameters and the influence of the uncertainty is estimated.

When the model contains multiple scales and high-contrast information, the computation of the snapshots may be quite expensive since we are required to resolve all scales of the model in a very fine mesh. To this end, we employ mixed GMsFEM  to compute snapshots in a coarse grid and develop reduced mixed GMsFE basis methods. To get a set of optimal reduced multiscale basis functions from the snapshots, we propose two sampling strategies: basis-oriented cross-validation (BOCV) and proper orthogonal decomposition (POD). For the two approaches, we select some optimal samples from a training set by a  greedy algorithm  \cite{gerner2012certified} based on mixed GMsFEM.  Then we
 construct the snapshots for the reduced mixed GMsFE basis functions  based on the optimal sample set. The optimal reduced mixed GMsFEM basis functions are identified by an incremental constructive manner from the snapshots using BOCV and POD. In the BOCV method, the optimal basis functions are searched with a minimum average error for a validation set. The POD is devoted to finding a low rank approximation to the space spanned by snapshots, and
 has been widely used in model reduction.  We use the two different approaches to generate a set of optimal reduced basis functions, and construct a reduced order model by projecting the full order model to the space spanned by the optimal reduced basis functions. We carefully compare BOCV with POD and find that BOCV can achieve better accuracy and robustness than POD in the reduced mixed GMsFE basis methods.

Although the reduced mixed GMsFE basis model  needs much less computation effort than original full order model, it may be not a very small-scale problem because the reduced order model still involves many unknowns.
This is still not desirable for the online computation in a many-query situation.
 In order to significantly improve the online computation, we want to get a reduced model representation for the solution and model outputs, which can accurately express the reduced model and can be used to estimate the influence of the uncertainty directly. To get the representation, we propose two methods:  least-squares method of snapshots (LSMOS) and sparse tensor approximation  (STA). They attempt to construct
 the representation with a form of variable-separation, which has been explored  in many applications \cite{beylkin2005algorithms, gonzalez2010recent, doostan2009least, doostan2007least}.
 STA is particularly useful  when the model
 suffers from the difficulty of high-dimensionality.

\begin{figure}[htbp]
       \centering
       \includegraphics[width=5in, height=3in]{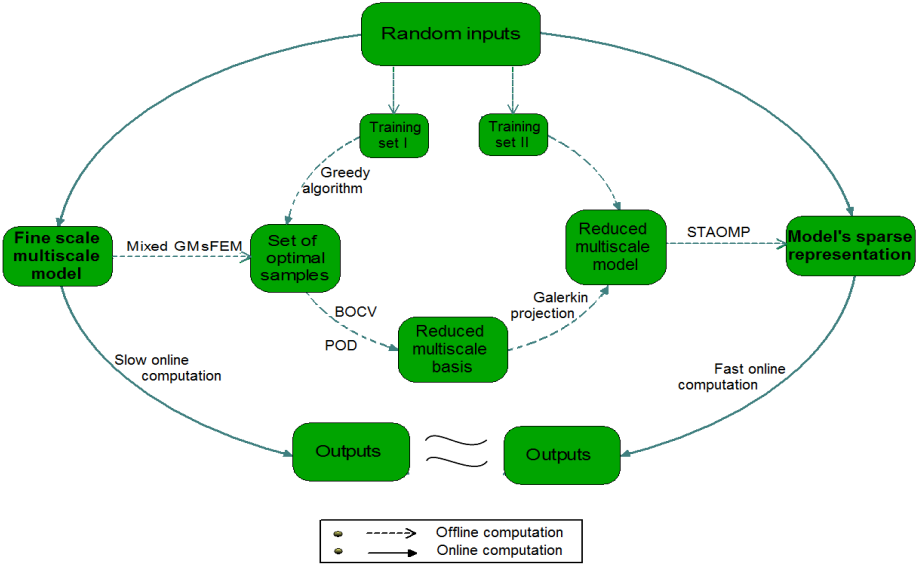}
       \caption{Schema of model's sparse representation based on reduced mixed GMsFE basis methods}
      \label{Schema}
\end{figure}

 LSMOS is based on Karhunen-Lo\`{e}ve expansion (KLE) \cite{loeve1978probability} and can extract most important information from a set of snapshots.  When it is difficult to get the knowledge of a covariance function for the
 model's ouputs, we can use method of snapshots \cite{handler2006karhunen} to construct a KL expansion for the model's outputs.
  In LSMOS, the coefficient functions of random variables in the KLE are the functions of the random parameters of inputs.
  We can use orthogonal polynomial basis functions to express the coefficient functions of the random variables.  To determine the coordinates of the orthogonal polynomial basis functions,
  we  take a few snapshots and use least-squares methods.  If the dimension of the random parameters is high,  a large number of snapshots are required for the least-squares methods.  Thus,
  the computation may be expensive.  To overcome the difficulty,
  we here propose a sparse tensor approximation method to obtain a sparse representation of the model's outputs instead of using   LSMOS.
  The sparse tensor approximation exploits the inherent sparsity of the model and the reduced GMsFE basis.
  The number of the  effective basis functions for many practical  models  is often  small. We use the optimization methods from compressive sensing \cite{rish2014sparse}  to extract  the effective basis functions and obtain
  a sparse representation.  There are roughly  two classes  of approaches  to obtain the sparse representation: optimization based on $l_0$-norm and  convex optimization  \cite{tropp2007signal, elad2009sparse, rish2014sparse}. The typical  methods of the $l_0$ optimization include  orthogonal matching pursuit (OMP) and iterative hard thresholding. The convex optimization based on $l_1$-norm includes least angle regression, coordinate descent and proximal methods.  In this work, we get the sparse tensor approximation using orthogonal matching pursuit, i.e., STAOMP.
  We integrate the reduced mixed GMsFE basis methods and STAOMP together to get the model's sparse representation in coarse scale.
  The proposed  model's sparse representation method can  significantly enhance the computation efficiency for the multiscale models
  with high-dimensional random inputs. The Figure \ref{Schema} illustrates the procedure to get the model's sparse representation based on reduced mixed GMsFE basis methods.

The paper is organized  as follows. In the next section, we give some preliminaries and notations for the paper. Section \ref{ssec:PGDDF-LSMOS-STAOMP} is devoted to developing LSMOS and STAOPM to represent a generic random field.
In Section \ref{ssec:Re-MIXGMsFEM-RB-method}, we  present the construction of the mixed GMsFE basis functions and introduce the reduced mixed GMsFE basis methods.  Section \ref{ssec:sampling algorithm} is to present two approaches  for constructing the optimal reduced basis: BOCV and POD. In Section \ref{ssec:Online-online}, we combine  STAOMP with the reduced mixed GMsFE basis methods together and get the model's sparse representation. A few  numerical examples are provided  in Section \ref{ssec:Numerical result} to illustrate the performance of all the methods developed above. In the last section, we make some comments and conclusions.
%%%%%%%%%%%%%%%%%%%%%%%%%%%%%%%%%%%%%%%%%%%%%%%

\section{Preliminaries and notations}

%\label{parametric PDEs}
\label{ssec:prelim}
In this section, we present some preliminaries and notations for the rest of paper. Let $L^2(D)$ be the space of square integrable functions over a domain $D$ with the $L^2$-norm $\|\cdot\|_{L^2(D)}^2=(\cdot, \cdot)$.
We define the space
\[
V:=H(\text{div},D)=\big\{u\in [L^2(D)]^d: \text{div}(u)\in L^2(D)\big\}.
\]
The space $V$ is a Hilbert space with the norm given by $\|u\|_{H(\text{div},D)}^2=\|u\|_{L^2(D)}^2+\|\text{div}(u)\|_{L^2(D)}^2$.
We consider parameterized elliptic PDEs in mixed formulation (\ref{mix-eq})
with the Neumann boundary condition $k\nabla p\cdot n=g$ on $\partial D$, where $n$ is the outward unit-normal vector on $\partial D$.
For a simplicity of notation, we denote $L^2(D)$ by $Q$ and
\[
V^0=H_0(\text{div},D):=\big\{u\in H(\text{div},D): u\cdot n=0 ~\text{on}~ \partial D\big\}.
\]
The problem (\ref{mix-eq}) leads to  the variational problem: for any parameter $\mu\in \Gamma\subset \mathbb{R}^p$, we find $\{v(\mu),p(\mu)\}\in V \times Q$ such that $v(\mu) \cdot n =g$ on $\partial D$,
\begin{eqnarray}
\label{mix-eq-weak}
\begin{cases}
\begin{split}
a\big(v(\mu),u; \mu\big)-b\big(u, p(\mu); \mu\big)&=0 ~~~\forall ~u \in V^0, \\
b\big( v(\mu), q; \mu\big)&=l(q) ~~~\forall ~q \in Q.
\end{split}
\end{cases}
\end{eqnarray}
where
\[
a(v,w):=\int_{D}k(x,\mu)^{-1}(v\cdot w )dx, ~~b(v,w):=\int_{D}\nabla\cdot(v)qdx,~~\text{and}~~ l(q):=\int_{D}f qdx.
\]
For any parameter $\mu\in \Gamma$,  $a(\cdot,\cdot;\mu):V\times V\longrightarrow \mathbb{R}$ is a symmetric bilinear form , $b(\cdot,\cdot;\mu):V \times Q\longrightarrow \mathbb{R}$ is bilinear form, and $l(\cdot;\mu)$ is a bounded linear functional over $Q$.
Given $\mu\in \Gamma$, we evaluate the output $G(x,\mu)$ of the model (\ref{mix-eq-weak}),
\begin{eqnarray*}
\label{evaluate-eq}
G(x,\mu)=L\big(v(x,\mu)\big),
\end{eqnarray*}
where $L$  is a operator  on $V$.

For well-posedness of (\ref{mix-eq-weak}), we assume that $a(\cdot , \cdot; \mu)$ is continuous and coercive over $V$ for all $\mu \in \Gamma$, i.e., there exist $\gamma>0$ and $\alpha >0$ such that
\begin{eqnarray}
\label{bounded-a}
\gamma(\mu):=\sup_{u \in V}\sup_{v \in V}\frac{a(u,v;\mu)}{\|u\|_V \|v\|_V}<\infty,~~~ \forall \mu\in \Gamma,
\end{eqnarray}
\begin{eqnarray}
\label{coercive-a}
\alpha(\mu):=\inf_{u \in V}\frac{a(u,u;\mu)}{\|u\|_V^2}> 0,~~~ \forall \mu\in \Gamma,
\end{eqnarray}
and $b(\cdot , \cdot; \mu)$ is continuous
\begin{eqnarray*}
\label{bounded-b}
b(u,q)\leq\gamma_b\|u\|_V \|q\|_Q,   \quad \forall u \in V,~ q \in Q, ~~\forall \mu\in \Gamma,
\end{eqnarray*}
and satisfies the inf-sup condition
\begin{eqnarray}
\label{inf-sup condition}
 \beta(\mu):=\inf_{q\in Q}\sup_{u\in V}\frac{b(u,q ;\mu)}{\|u\|_V\|q\|_Q}> 0,   ~\forall \mu\in \Gamma.
\end{eqnarray}
%%By (\ref{bounded-a}) and (\ref{coercive-a}), the bilinear form $a(\cdot , \cdot; \mu)$ provides a norm $\|\cdot\|_{E}^2=a(\cdot,\cdot;\mu)$ on $V$.

To fulfill  offline-online computation,  we assume that the parametric bilinear form $a(\cdot, \cdot; \mu)$ is affine with respect to $\mu$, i.e.,
\begin{eqnarray}
\label{affinely-ag}
a(u,v;\mu)&=\sum_{i=1}^{m_{a}}k^{i}(\mu)a^{i}(u,v), \quad \forall u,v\in V,\quad \forall \mu \in \Gamma,
\end{eqnarray}
where $k^{i}: \Gamma\longrightarrow \mathbb{R}$ is a $\mu$-dependent function and $a^{i}: V \times V\longrightarrow \mathbb{R}$ is a symmetric bilinear form independent of $\mu$, for each $i= 1,\cdots , m_{a}$. The affine assumption (\ref{affinely-ag}) is crucial to achieve the decomposition of  offline-online computation   for  many queries to model's outputs. When $a(\cdot, \cdot; \mu)$ is not affine with regard to $\mu$, such an expansion can be obtained by using some variable separation methods for $a(\cdot, \cdot; \mu)$, e.g.,  EIM \cite{eftang2010posteriori}, LSMOS and  STAOMP. The LSMOS and  STAOMP will be introduced  in Section \ref{ssec:PGDDF-LSMOS-STAOMP}.

Let $V_h\times Q_h$ be the pair of standard Raviart-Thomas space for the approximation of (\ref{mix-eq-weak}) on a fine grid $\mathcal{K}^h$.
Then the standard mixed FEM  of problem (\ref{mix-eq-weak}) reads: given any $\mu \in \Gamma$, we find $\{v_h(\mu),p_h(\mu)\}\in V_h \times Q_h$ such that
\begin{eqnarray}
\label{mix-FEM-dis}
\begin{cases}
\begin{split}
a(v_h(\mu),u; \mu)-b(u, p_h(\mu); \mu)&=0 ~~~\forall ~u \in V_h^0, \\
b( v_h(\mu), q; \mu)&=l(q) ~~~\forall ~q \in Q_h,
\end{split}
\end{cases}
\end{eqnarray}
where  $V_h^0=V_h\bigcap\{v\in V_h: v\cdot n=0 ~\text{on}~\partial D \}$.
We remark that the fine-grid solution is considered as a reference solution in the paper.

\section{ Variable-separation  strategies}
\label{ssec:PGDDF-LSMOS-STAOMP}

In order to achieve offline-online computation, we want to represent model's  inputs (e.g., coefficients and source terms) and outputs
by a form of variable-separation.   Let $G(x,\mu)$ be a generic function parameterized by $\mu$.
In this section, we introduce two strategies, i.e., LSMOS and STAOMP, to construct an approximation in the form
\begin{eqnarray}
\label{separating presentation}
G(x,\mu)\approx G_N(x,\mu):=\sum_{i=1}^{N}\zeta_i(\mu)g_i(x),
\end{eqnarray}
where  $\zeta_i(\mu)$ only depends on $\mu$ and $g_i(x)$ only depends on $x$.
Let $L_2(\Gamma)$ be the space of square integrable functions over the parameter space $\Gamma$. We consider the approximation in  the tensor space $\mathcal{H}\otimes L_2(\Gamma)$, where
$\mathcal{H}$ is a Hilbert space with regard to the spatial space.
We define an inner product on $\mathcal{H}\otimes L_2(\Gamma)$ by
\begin{eqnarray*}
\label{TENSOR PRODUCT}
(w,u)_{\mathcal{H}\otimes L_2(\Gamma)}=E[(w,u)_\mathcal{H}]:=\int_{\Gamma}(w,u)_\mathcal{H} \rho(\mu)d\mu,
\end{eqnarray*}
where $\rho(\mu)$ is the density function for the random parameter $\mu$. The norm is defined by
\[
\|u\|_{\mathcal{H}\otimes L_2(\Gamma)}^2:=(u,u)_{\mathcal{H}\otimes L_2(\Gamma)}.
\]
We assume that $S_N \subset \mathcal{H}\otimes L_2(\Gamma)$ is a finite dimensional subspace  space, and $\{\Psi_i\}_{i=1}^N$ is a set of basis functions for $S_N $.
We want to find an approximation of $G(x, \mu)$ in $S_N$ such that
\begin{eqnarray}
\label{output-approximation}
 \|G(x,\mu)-\sum_{j=1}^{N} c_{j}\Psi_j(x,\mu) \|_{\mathcal{H}\otimes L_2(\Gamma)}\leq \delta,
\end{eqnarray}
where $\delta$ is a  given threshold.

\subsection{Least-squares method of snapshots}

In this subsection, we introduce a least-squares method of snapshots to get the approximation (\ref{output-approximation}).
Let $\Xi_{t}$ be a collection of a finite number of samples in $\Gamma$ and the cardinality $|\Xi_{t}|=n_t$.
For $\forall ~\mu\in \Xi_{t}$, we can split $G(x,\mu)$ into two parts, i.e.,
\[
G(x,\mu)=\bar{G}(x)+\tilde{G}(x,\mu),
\]
where $\bar{G}(x):=E[G(x,\cdot)]=\frac{1}{n_{t}}\sum_{i=1}^{n_{t}}G(x,\mu_i)$ is the mean, and $\tilde{G}(x,\mu)=G(x,\mu)-\bar{G}(x)$ is a random fluctuating part.
To obtain $\tilde{G}(x,\mu)$, we take a set of snapshots $\{\tilde{G}(x,\mu_i)\}_{i=1}^{n_{t}}$ and compute a covariance matrixes $\textbf{C}$, whose entries can be defined by
\[
\textbf{C}_{n,m}:=\frac{1}{n_{t}}\bigg(\tilde{G}(x,\mu_n),\tilde{G}(x,\mu_m)\bigg)_\mathcal{H}.
\]
Let $\{\hat{\lambda}_k,  \textbf{e}_k\}$  be the eigen-pairs (normalized)  of $\textbf{C}$, $1\leq k\leq n_{t}$. Set $(\textbf{e}_k)_j=e_k^{j}$, we define the functions
\begin{eqnarray}
\label{KLE basis}
g_k(x):=\frac{1}{\sqrt{\hat{\lambda}_k n_{t}}}\sum_{j=1}^{n_{t}}e_k^{j}\tilde{G}(x,\mu_j).
\end{eqnarray}
It is easy to get $(g_k,g_l)_\mathcal{H}=\delta_{k,l}$, $1\leq k$, $l\leq n_{t}$. Then it holds that
\[
\tilde{G}(x,\mu)\approx \sum_{i=1}^{ n_{t}} \sqrt{\hat{\lambda}_i} \zeta_i(\mu) g_i(x),
\]
where $\{\zeta_i(\mu)\}_{i=1}^{n_{t}}$  are given by
%are uncorrelated random variables with zero mean and unit variance, and
\begin{eqnarray}
\label{KLE parameter}
\begin{split}
\zeta_i(\mu):= \frac{1}{\sqrt{\hat{\lambda}_i}} \big(\tilde{G}(\cdot ,\mu),g_i\big)_\mathcal{H}.
\end{split}
\end{eqnarray}
Thus we get the decomposition
\begin{eqnarray}
\label{KLE Snapshots}
G(x,\mu)\approx \bar{G}(x)+\sum_{i=1}^{ M} \sqrt{\hat{\lambda}_i} \zeta_i(\mu) g_i(x).
\end{eqnarray}

 Because the equation (\ref{KLE parameter}) involves  $\tilde{G}(\cdot,\mu)$, it can not be directly used to compute the functions $\{\zeta_i(\mu)\}_{i=1}^{n_{t}}$ for arbitrary  $\mu$.
  We will apply least-squares methods    to approximate  $\{\zeta_i(\mu)\}_{i=1}^{n_{t}}$ based on  orthogonal polynomials.
  Let $\{p_i(\mu)\}_{i=1}^{M_{g}}$  be the set of orthogonal polynomials basis functions with total degree less than $N_g$ regarding to the parameter variable $\mu$.   We rearrange the set of basis functions  from the first to the last one and place them  in the following row vector,
 \[
 \big[p_1(\mu), p_2(\mu),...,p_{M_{g}}(\mu)\big], ~~~~\text{where} ~~M_{g}=\left( \begin{array}{cc}  N_g+p\\  N_g \\\end{array}\right).
 \]
 For the sample date  $\Xi_{t}$, we compute $[p_1(\mu_j), p_2(\mu_j),\cdots,p_{M_{g}}(\mu_j)]$ and $\zeta_i(\mu_j)=\frac{1}{\sqrt{\hat{\lambda}_i}} (\tilde{G}(\cdot,\mu_j),g_i)_\mathcal{H}$ ($j = 1,\cdots,n_{t}$).  They are putted in the following matrix $\textbf{A}$ and vector $\textbf{F}$ , respectively,
 \begin{eqnarray}
\label{least square A}
 \textbf{A}:=\left[ \begin{array}{ccc}  p_1(\mu_1)& \ldots & p_{M_{g}}(\mu_1)\\ \vdots &\ddots &\vdots\\ p_1(\mu_{n_{t}})& \ldots & p_{M_{g}}(\mu_{n_{t}})\\\end{array}\right] ,
 \end{eqnarray}
 \begin{eqnarray}
\label{least square F}
 \textbf{F}:=[\zeta_i(\mu_1) \cdots  \zeta_i(\mu_{n_{t}})]^{T}.
 \end{eqnarray}
 We obtain the approximation of the parameter functions $\zeta_i(\mu)$ by solving the following least square problem,
\begin{eqnarray}
\label{least square problem}
\begin{split}
\textbf{d}=\arg\min_{\alpha}\|\textbf{A}\alpha-\textbf{F}\|_2.
\end{split}
\end{eqnarray}
 Thus we get $\zeta_i(\mu)\approx \sum_{i=1}^{M_{g}}d_ip_i(\mu)$, and $d_i=(\textbf{d})_i$.

\subsection{ Sparse tensor approximation}

In this subsection, we are devote to seeking the optimal solution $\textbf{c}$ for (\ref{output-approximation}) with the minimum number of non-zero terms. This can be formulated as the optimization problem:
\begin{eqnarray}
\label{P0-PROBLEM}
\arg\min_{\textbf{c}}\|\textbf{c}\|_0 ~~ \text{subject to} ~~\|G(x,\mu)-\sum_{j=1}^{N} c_{j}\Psi_j(x,\mu) \|_{\mathcal{H}\otimes L_2(\Gamma)}\leq \delta,
\end{eqnarray}
where
\[
\textbf{c}=(c_1, \cdots, c_N), \quad \text{and} \quad         \|\textbf{c}\|_0=\sharp\{j:c_{j}\neq 0\}.
\]
Thus the sparse solution  can be constructed as follows.\\
  $\bullet$ \textit{Step 1: Find the optimal $N$-dimensional subspace $\mathcal{H}_N\subseteq \mathcal{H}$}

With the snapshots $\{G(x,\mu_i)\}_{i=1}^{n_{t}}$, we can use methods of snapshots  or POD   to construct the optimal $N$-dimensional subspace
\[
\mathcal{H}_N=\text{span}\big\{g_j(x): 1\leq j\leq N\big\}.
\]
$\bullet$ \textit{Step 2: Choose orthogonal polynomials for $L_2(\Gamma)$ and construct the finite dimensional approximation space $S_{N\times M}\subseteq \mathcal{H}\otimes L_2(\Gamma)$}

For the parameter, we introduce appoximation spaces $X_M$ for $L_2(\Gamma)$,
\[
X_M=\text{span}\big\{p_i(\mu): 1\leq i\leq M\big\}.
\]
A finite dimensional approximation space $S_{N\times M}\subseteq \mathcal{H}\otimes L_2(\Gamma)$ is then obtained by
\[
S_{N\times M}:=\mathcal{H}_N\otimes X_M=\text{span}\big\{p_i(\mu)g_j(x): 1\leq i\leq M,1\leq j \leq N\big\}.
\]
To simplify notation, we use the following single-index notation
\[
S_{N\times M}=\big\{w(x,\mu)=\sum_{i\in I}w_i\Psi_i(x,\mu); ~w_i\in \mathbb{R}\big\},
\]
where $I=\{1,\cdots,N\}\times \{1,\cdots,M\}$ and $\Psi_i(x,\mu)=p_{i_1}(\mu)u_{i_2}(x)$.\\
$\bullet$ \textit{Step 3: Construct the sparse solution based on Orthogonal-Matching-Pursuit}

In general, the optimization problem (\ref{P0-PROBLEM}) is an NP-hard problem. We attempt to seek efficient algorithms to approximately solve (\ref{P0-PROBLEM}).
 There are a few approaches \cite{rish2014sparse} to solve the problem (\ref{P0-PROBLEM}).   In the paper, we focus on  Orthogonal Matching Pursuit (OMP) algorithm \cite{tropp2007signal} to get a sparse solution.

We place the basis functions $\{\Psi_i(x,\mu)\}_{i=1}^{M\times N} \subseteq S_{N\times M}$ in the following row vector,
\[
\Psi:=[\Psi_1, \Psi_2,\cdots,\Psi_{M\times N}].
\]
Assume that $n=n_x\times n_{\mu}$ sample data $\{(x_i,\mu_j)\}\subseteq D \times \Gamma$ ($1\leq i\leq n_x$, $1\leq j\leq n_{\mu}$) are chosen to solve the optimization problem (\ref{P0-PROBLEM}). For the sample date, we compute $\Psi(x_i,\mu_j)$ and $G(x_i,\mu_j)$, and they
 are putted in the following matrix $\Pi\in\mathbb{R}^{n\times (M\times N)}$ and vector $\textbf{b}\in\mathbb{R}^{n\times 1}$, respectively,
\begin{eqnarray}
\label{Pi-b-matrix}
 \Pi:=\left[ \begin{array}{ccc}  \Psi_1(x_1,\mu_1)& \ldots & \Psi_{M\times N}(x_1,\mu_1)\\\vdots &\ddots &\vdots\\ \Psi_1(x_{n_x},\mu_{n_{\mu}})& \ldots &\Psi_{M\times N}(x_{n_x},\mu_{n_{\mu}})\\\end{array}\right], ~~~~\textbf{b}:=\left[ \begin{array}{ccc}  G(x_1,\mu_1)\\\vdots\\ G(x_{n_x},\mu_{n_{\mu}}) \\\end{array}\right].
\end{eqnarray}
 The coefficient vector can be solved by the following optimization problem,
\begin{eqnarray}
\label{P0-PROBLEM-discrete}
\arg\min_{\textbf{c}}\|\textbf{c}\|_0 ~~ \text{subject to} ~~\|\textbf{b}-\Pi \textbf{c}\|_2\leq \delta.
\end{eqnarray}
The main idea of OMP is to pick columns in a greedy manner. At each iteration, we choose the column of $\Pi$ that is most strongly correlated with the remaining part of $\textbf{b}$, i.e., the residual $\textbf{r}_k$ in Algorithm \ref{algorithm-OMP}. Then we subtract off the contribution to $\textbf{b}$ and iterate with regard to the updated residual. After $k$ iterations, the algorithm can identify  the correct set of columns.
%Let $ \emph{I}$ be the set of the indexes corresponding to the
The OMP is described in Algorithm \ref{algorithm-OMP}.
\begin{algorithm}
\caption{ Orthogonal-Matching-Pursuit}
     \textbf{Input}: A matrix $\Pi$, the vector $\textbf{b}$ and the error tolerance $\varepsilon$\\
     \textbf{Output}:  The sparse solution $\textbf{c}$ and the solution support $ \emph{I}:=supp(\textbf{c})$ \\
      ~1:~~Initialize the residual $\textbf{r}_0=\textbf{b}$, the index set $\emph{I}^0={\O}$, the iteration counter $k=1$, and \\
      $~~~~~\Pi_0$ is an empty matrix;\\
      ~2:~~Find the index $j_0$ that solves the easy optimization problem:\\
      $~~~~~j_0=\arg\max_{j=1,...,M\times N}|\langle \textbf{r}_k, \Psi_j\rangle|$;\\
      ~3:~~Update the index set $\emph{I}^k=\emph{I}^{k-1}\cup \{j_0\}$, the matrix $\Pi_k=[\Pi_{k-1} ~\Psi_{j_0}]$; \\
      ~4:~~Solve a least-squares problem to obtain a new estimate: \\
      $~~~~~\textbf{c}^k=\arg\min_{\textbf{c}}\|\textbf{b}-\Pi_k \textbf{c}\|_2$;\\
      ~5:~~Calculate the new residual: $\textbf{r}_k=\textbf{b}-\Pi_k \textbf{c}^k$; \\
      ~6:~~$k\rightarrow k+1$, return to Step 2 if $\frac{\|\textbf{r}_k\|_2}{\|\textbf{b}_k\|_2}\geq \varepsilon$, otherwise \textbf{terminate} .\\
      ~7:~~$\text{Mt}=k$, $\emph{I}=\emph{I}^k$, and $\textbf{c}=\textbf{c}^k$.
      \label{algorithm-OMP}
\end{algorithm}
The residual is always orthogonal to the columns that have been selected. In fact, we can get the conclusion from step $4$ in Algorithm \ref {algorithm-OMP} by vanishing the derivative of $\|\textbf{b}-\Pi_k \textbf{c}\|_2$, i.e.,
\[
-\Pi_k^{T}(\textbf{b}-\Pi_k \textbf{c})=-\Pi_k^{T}\textbf{r}_k=0.
\]
Thus, OMP never selects the same column twice. Provided that the residual is nonzero, the algorithm selects a new atom at each iteration and the matrix $\Pi_k$ has full column rank.

\begin{rem} For the classical least-squares method, it is required that the number of parameter sample scales quadratically with the number of unknowns \cite{chevreuil2015least}. However, STAOMP can provide an accurate  approximation by a much fewer number of samples.
\end{rem}

%%%%%%%%%%%%%%%%%%%%%%%%%%%%%%%%%%%%%%

\section{Reduced mixed GMsFE basis method}
\label{ssec:Re-MIXGMsFEM-RB-method}

Let $\mathcal{K}^H$ be a  conforming coarse partition for the computational domain $D$, where $H$ is the coarse  mesh size. Each coarse-grid block is further partitioned into a connected union of fine-grid blocks, we get the fine grid partition $\mathcal{K}^h$. Let $\varepsilon^H:=\bigcup_{i=1}^{N_e}\{E_i\}$  be the set of all edges/interfaces of coarse mesh $\mathcal{K}^H$ and $N_e$ the number of coarse edges. The coarse neighborhood $w_i$ corresponding to the coarse edge $E_i$ is defined by
\[
w_i=\bigcup\{K_j\in\mathcal{K}^H; E_i\in\partial K_j\}.
 \]
 In the mixed GMsFEM, the velocity field is approximated by using mixed GMsFE basis functions, while piecewise constant functions over  $\mathcal{K}^H$ are used to approximate the pressure field. Let $Q_H$ be the space of piecewise constant functions over the coarse grid $\mathcal{K}^H$. Let $\{\phi_j\}_{j=1}^{L_i}$ be the set of GMsFE basis functions corresponding to edge $E_i$. We define the GMsFE space for the velocity field as the linear span of all local basis functions, i.e.,
 \[
 V_H=\bigoplus_{\varepsilon^H}\{\phi_j\}_{j=1}^{L_i}.
 \]
 Let $V_H^0=V_H\bigcap\{u\in V_H: u\cdot n=0 ~\text{on}~\partial D \}$ be a subspace of $V_H$. Thus the mixed GMsFEM is to find $\{v_H(\mu), p_H(\mu)\}\in V_H\times Q_H$ such that
\begin{eqnarray}
\label{discre-mix-eq}
\begin{cases}
\begin{split}
a(v_H(\mu),u; \mu)-b(u, p_H(\mu); \mu)&=0 ~~~\forall~ u \in V_H^0, \\
b(v_H(\mu), q; \mu)&=l(q) ~~~\forall ~q \in Q_H,\\
\end{split}
\end{cases}
\end{eqnarray}
where $v_H(\mu)\cdot n=g_H ~\text{on} ~ \partial D$, and for each coarse edge $E_i \in \partial D$,
\[
\int_{E_i} (g_H-g) \phi_j \cdot n=0, ~ j=1,\cdots, L_i.
\]
It is easy to see that $Q_H\subset Q_h$ and $V_H\subset V_h$ in the mixed GMsFEM. We will briefly present the mixed GMsFEM and
introduce  reduced mixed GMsFE basis methods.

\subsection{Mixed GMsFE space}
In this section, we follow mixed GMsFEM \cite{cel15} and  present the construction of the GMsFE space $V_H$ for the approximation of the velocity field. We first t generate the snapshot space and
 then use spectral decomposition  to obtain a lower-dimensional offline space.
%$\mathbf{Snapshot ~~space}$\\

\begin{figure}[htbp]
\centering
  \includegraphics[width=3.5in, height=3in]{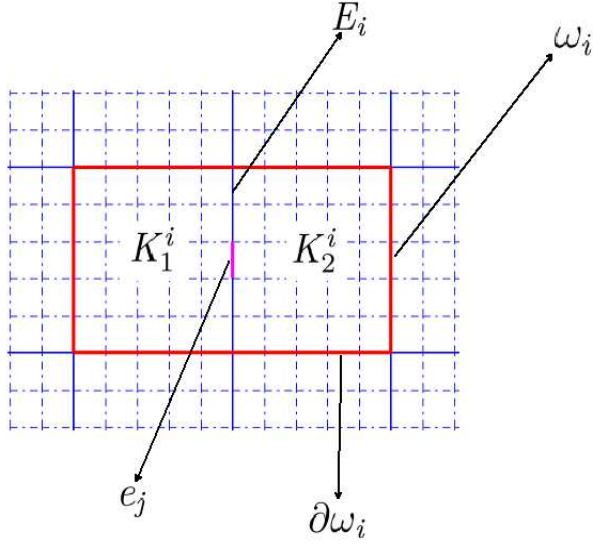}
  \caption{Illustration of  a neighborhood $\omega_i = K_1^i\cup K_2^i$.}
  \label{fig1-coarseblock}
\end{figure}

Let $E_i\in \varepsilon^H$ be a coarse edge and $e_j\subset E_i$ a fine edge, and define a piecewise constant function $\delta_j^i$ on $E_i$ as
\begin{eqnarray*}
\delta_{j}^{i}&=
\begin{cases}
\begin{split}
1,& ~\text{on} ~ e_{j}, \\
0,&~ \text{on other fine  edges of}~ E_i.
\end{split}
\end{cases}
\end{eqnarray*}
We solve the following problem on the coarse neighborhood  $w_i$ corresponding to the edge $E_i$,
\begin{eqnarray}
\label{local-mix-eq-snap}
\begin{cases}
\begin{split}
k(x,\mu)^{-1}v_j^i(\mu)+\nabla p_j^i(\mu)&=0 ~~\text{in} ~~w_i, \\
\nabla\cdot v_j^i(\mu)&=\alpha_j^i ~\text{in} ~~w_i, \\
 v_j^i(\mu)\cdot n_i&=0  ~~\text{on}  ~~\partial w_i,
\end{split}
\end{cases}
\end{eqnarray}
where the constant $\alpha_j^i$ satisfies the compatibility condition $\int_{K_n}\alpha_j^i=\int_{E_i}\delta_j^i$ for all $K_n\subset w_i$, and $n_i$ denotes the outward unit normal vector on $\partial w_i$.  The local problem (\ref{local-mix-eq-snap}) is solved separately in the coarse-grid blocks of $w_i$. We need an extra boundary condition on $E_i$ for well-posedness. Let $J_i$ be the total number of fine-grid edges on $E_i$ and $E_i=\cup_{j=1}^{J_i}e_j$, where  $e_j$ denotes a fine-grid edge.
 The remaining boundary condition on the coarse edge $E_i$ for the local problem (\ref{local-mix-eq-snap}) is taken as
\begin{eqnarray*}
\label{extra boundary}
v_j^i(\mu)\cdot m_i=\delta_j^i  ~~\text{on}  ~~E_i ,
\end{eqnarray*}
where $m_i$ is a fixed unit-normal vector on $E_i$. See Figure \ref{fig1-coarseblock} for illustration of a coarse neighborhood.
Then we define the snapshot space $V_{snap}$ space by
\[
V_{snap}=\text{span}\{v_j^i(\mu):1\leq j\leq J_i,1\leq i\leq N_e\}.
\]

%$\mathbf{Offline ~~space}$\\
After the construction of the snapshot space, the offline space is constructed by performing some local spectral problem on the snapshot space.
Let $V_{snap}^i$ be the snapshot space corresponding to the coarse edge $E_i$, i.e.,
\[
V_{snap}^i=\text{span}\{v_j^i(\mu):1\leq j\leq J_i\}.
\]
We consider the local spectral problem: find eigenpair  $\{\lambda, v\}$ ($v\in V_{snap}^i$) such that
\begin{eqnarray}
\label{tem-spectral-problem}
a_i(v,w)=\lambda s_i(v,w)~~~~\forall~ w \in V_{snap}^i,
\end{eqnarray}
where
\begin{eqnarray*}
\label{bilinear-form}
\begin{split}
a_i(v,w)&=\int_{E_i}k(x,\mu)^{-1}(v\cdot m_i)(w\cdot m_i)dx,\\
s_i(v,w)&=\int_{w_i}k(x,\mu)^{-1}v\cdot wdx+\int_{w_i}(\nabla\cdot v) (\nabla\cdot w)dx.
\end{split}
\end{eqnarray*}
Suppose that the eigenvalues of (\ref{tem-spectral-problem}) are arranged in increasing order, the corresponding eigenvectors are denoted by $z_r^i$.
We take the eigenfunctions $\{\phi_r^i\}_{r=1}^{l_i(\mu)}$ corresponding to the first $l_i(\mu)$ eigenvalues to form the local reduced snapshot space,
\[
V_{\text{off}}^i=\text{span}\{\phi_j^i(\mu):1\leq j\leq l_i(\mu)\}.
\]
Then the offline space is
\[
V_{\text{off}}=\text{span}\{\phi_j^i(\mu):1\leq j\leq l_i(\mu),1\leq i\leq N_e\}.
\]
To simplify notation, we use the following single-index notation
\begin{eqnarray*}
\label{local-snapshot-space}
V_{\text{off}}=\text{span}\{\phi_k(\mu):1\leq k \leq n_{\text{off}}\},
\end{eqnarray*}
where $n_{\text{off}}=\sum_{i=1}^{N_e}l_i(\mu)$.

\subsection{Reduced mixed GMsFE space }

In the reduced basis method, we construct a set of reduced basis for arbitrary parameter $\mu$  based on a set of snapshots.
 If the snapshots have strong multiscale features, we have to use a very fine mesh to resolve the features in all scales. This computation may be very expensive. To overcome the difficulty, we use mixed GMsFEM to compute the snapshots.
Let $\Xi_{op}$ be an optimal parameter set, which is a collection of a finite number of samples in $\Gamma$. In the paper,  we will use a greedy algorithm to identify $\Xi_{op}$, which will be discussed in Section \ref{greedy selection}.

The reduced mixed GMsFEM is devoted to approximating the solution $\{v(\mu), p(\mu)\}$  of the problem (\ref{mix-eq-weak}) by a set of pre-computed basis functions $\{\psi_j^i:1\leq j\leq N,1\leq i\leq N_e\}$, which are selected from
\[
\Sigma:=\{\phi_j^i(\mu):\mu\in\Xi_{op}, 1\leq j\leq l_i,1\leq i \leq N_e\}
 \]
 with some optimal strategies. Let $V_H^N$ be an $(N \times N_e)$-dimensional subspace of $V$ and define
\[
\mathcal{F}:=\{v(\mu)\in V_h:\mu\in \Gamma\}.
 \]
To assess approximation property, it is natural to compare the subspace $V_H^N$  with the best $N \times N_e$-dimensional subspace spanned by some elements of $\Sigma$,  which minimizes the projection error for the  $\Sigma$ over all $N \times N_e$-dimensional subspaces of V. This minimal error can be described by the Kolmogorov width
\begin{eqnarray*}
\label{Kolmogorov width}
d_{N \times N_e}(Y_{N \times N_e}, \mathcal{F}):=\inf\{ E(\mathcal{F};Y_{N \times N_e}): Y_{N \times N_e} \text{ is an $N \times N_e$-dimensional subspace of }  V \},
\end{eqnarray*}
where $E(\mathcal{F};Y_{N \times N_e})$ is the angle between $\mathcal{F}$ and $Y_{N \times N_e}$ under a metric. We construct a finite dimensional space,  which is spanned by elements of $\Sigma$ with  good approximation. The procedure is described  as follows:\\

$\bullet$ $V_H^1=\mathop{\arg\min}\limits_{\mathop{Y_{1\times N_e}\subset V}\limits_{dim (Y_{1\times N_e)}=N_e}} d_{1\times N_e}(Y_{1\times N_e}, \mathcal{F})$,\\

$\bullet$ Assume that $V_H^{N-1}$ have been constructed. Then $V_H^{N}=\mathop{\arg\min}\limits_{\mathop{V_H^{N-1}\subset Y_{N \times N_e}\subset V}\limits_{dim (Y_{N \times N_e)}=N \times N_e}} d_{1\times N_e}(Y_{N \times N_e}, \mathcal{F})$.

A sequence of reduced GMsFE approximation spaces for velocity are obtained as follow,
\[
V_H^1\subset V_H^2\subset,. . .,\subset V_H^{N} \subset V_h,
\]
and a set of basis functions $\{\varphi_{n}^{i}: 1\leq i\leq N_e, 1\leq n\leq N\}$ are available.
 To obtain a set of $(\cdot,\cdot)_V$-orthonormal basis functions, we apply the Gram-Schmidt process to $\{\varphi_{n}^{i}: 1\leq i\leq N_e, 1\leq n\leq N\}$ in the $(\cdot,\cdot)_V$ inner product. The set of orthonormal basis functions is denoted by
 \begin{eqnarray}
\label{MsRB-basis}
\{\psi_{i}: 1\leq i\leq N \times N_e\}.
\end{eqnarray}
We note that the support of each basis function $\psi_{i}$  in (\ref{MsRB-basis}) is on a coarse block.

We apply Galerkin projection to construct a reduced model with the reduced multiscale basis functions defined in (\ref{MsRB-basis}).
Let  $\{v_H^N(\mu), p_H^N(\mu)\}\in V_H^N\times Q_H^N$   solve
\begin{eqnarray}
\label{RB-discre-mix-eq}
\begin{cases}
\begin{split}
a(v_H^N(\mu),u; \mu)-b(u, p_H^N(\mu); \mu)&=0 ~~~\forall~ u \in V_H^N, \\
b(v_H^N(\mu), q; \mu)&=l(q) ~~~\forall ~q \in Q_H,
\end{split}
\end{cases}
\end{eqnarray}
where $v_H^N(\mu)\cdot n=g_H ~\text{on} ~ \partial D$.

Suppose that $\{\psi_{i}\}_{i=1}^{N_e \times N}$ is the set of basis functions for $V_H^N$, and $\{\eta_{r}\}_{r=1}^{N_{\text{el}}}$ is the basis functions for $Q_H$, where $N_{\text{el}}$ is the number of the coarse element. Then solutions $\{v_H^N(\mu), p_H^N(\mu)\}$ can be represented by
\begin{eqnarray}
\label{GMsRB solution}
v_H^N(\mu)=\sum_{i=1}^{N_e \times N}v_{i}^N(\mu) \psi_{i}, ~~~p_H^N(\mu)=\sum_{r=1}^{N_{\text{el}}}p_r^N(\mu) \eta_{r}.
\end{eqnarray}
By plugging $u =\psi_{j}$, $1\leq j\leq N \times N_e$  and $q =\eta_{n}$, $1\leq n\leq N_{\text{el}}$ into (\ref{RB-discre-mix-eq}), we have
\begin{eqnarray}
\label{RB-stiffness-equations}
\begin{cases}
\begin{split}
\sum_{i=1}^{N_e \times N}a(\psi_{i},\psi_{j}; \mu)v_i^N(\mu)-\sum_{r=1}^{N_{\text{el}}} b(\psi_{j}, \eta_{r})p_r^N(\mu)&=0, ~~~\quad 1\leq j\leq N \times N_e,\\
\sum_{i=1}^{N_e \times N}b(\psi_{i}, \eta_{n})v_i^N(\mu)&=l(\eta_{n}),~\quad 1\leq n \leq N_{\text{el}}.
\end{split}
\end{cases}
\end{eqnarray}
Then the output of the  model can be evaluated by
\begin{eqnarray*}
\label{evaluate output}
G_H^N(\mu)&=L\big(v_H^N(\mu)\big).
\end{eqnarray*}

\subsection{Offline-online computation decomposition for reduced mixed GMsFE basis method}
\label{procedure}
The equation (\ref{RB-stiffness-equations}) implies a linear algebraic  system with $N \times N_e+N_{\text{el}}$ unknowns. It involves the computation of inner products with  entities $\{\psi_{i}\}_{i=1}^{N \times N_e}$ and $\{\eta_{n}\}_{n=1}^{N_{\text{el}}}$, each of which is represented by fine grid finite element basis functions of $V_h \times Q_h$.
This will lead to substantial computation for the input-output evaluation $\mu\longrightarrow G_H^N(\mu)$. With the assumption (\ref{affinely-ag}) of affine decomposition, the equation (\ref{RB-stiffness-equations}) can be rewritten by
\begin{eqnarray*}
\label{OffOnline matrix equations-GMsRB}
\begin{cases}
\begin{split}
\sum_{i=1}^{N_e \times N}\sum_{q=1}^{m_{a}}k^{q}(\mu)a^q(\psi_{i},\psi_{j})v_i^N(\mu)-\sum_{r=1}^{N_{\text{el}}} b(\psi_{j}, \eta_{r})p_r^N(\mu)=0, ~~~\quad 1\leq j\leq N \times N_e, \\
\sum_{i=1}^{N_e \times N}b(\psi_{i}, \eta_{n})v_i^N(\mu)=l(\eta_{n}),~\quad 1\leq n \leq N_{\text{el}}.
\end{split}
\end{cases}
\end{eqnarray*}
This gives rise to the matrix form
\begin{eqnarray}
\label{RB-G-matrix-form}
\begin{cases}
\begin{split}
\sum_{q=1}^{m_{a}}k^{q}(\mu)M_{N}^q\textbf{v}_N+B_N^T\textbf{p}_N=0,\\
B_N\textbf{v}_N=F_N,
\end{split}
\end{cases}
\end{eqnarray}
where
\[
(M_{N}^q)_{ij}=a^{q}(\psi_{i},\psi_{j}), \quad (B_N)_{in}=b(\psi_{i}, \eta_{n}), \quad (F_N)_{n}=l(\eta_{n}),\quad (\textbf{v}_N)_{i}=v_i^N,\quad (\textbf{p}_N)_{r}=p_r^N,
\]
\[
1\leq i,j\leq N \times N_e, ~~1\leq n\leq N_{\text{el}},~~1\leq r\leq N_{\text{el}}.
\]
 Because basis functions $\{\psi_{i}, \eta_{n}\}$ belong to the standard Raviart-Thomas space $V_h\times Q_h$ ,
 we can express them by
\[
\psi_{i}=\sum_{k=1}^{N_e^f}Z_{ik}\xi_{k}, 1\leq i\leq  N \times N_e, ~~ \eta_{n}=\sum_{k=1}^{N_{\text{el}}^f}I_{ik}\gamma_{k}, 1\leq n\leq N_{\text{el}},
\]
where $N_e^f$ is the number of fine edges and $N_{\text{el}}^f$ is the number of fine elements.
Let $(\mathcal{Z})_{ki}=Z_{ki}$ and $(\mathcal{I})_{ki}=I_{ki}$. Then we  get
\[
\mathbf{M}_{N}^{q}=\mathcal{Z}^{T}\mathcal{M}_{N_f}^{q}\mathcal{Z}, \quad  \mathbf{B}_{N}=\mathcal{I}^{T}\mathcal{B}_{N_f}\mathcal{Z}, \quad \mathbf{F}_{N}=\mathcal{I}^{T}\mathcal{F}_{N_f},
\]
where $(\mathcal{A}_{N_f}^{q})_{ij}=a^{q}(\xi_{j}, \xi_{i})$, $(\mathcal{B}_{N_f})_{ij}=b(\xi_{i},\gamma_{j})$, $(\mathcal{F}_{N_f}^{q})_{i}=l(\gamma_{i})$.
The matrixes $\mathcal{A}_{N_f}^{q},~\mathcal{B}_{N_f}$ and the vectors $\mathcal{F}_{N_f}^{q}$ are independent of parameter $\mu$,  and their  computation is once and  in  offline phase. The online computation
is to solve equation (\ref{OffOnline matrix equations-GMsRB}) and evaluate the output $G_H^N(\mu)$ for any $\mu\in \Gamma$.  The system (\ref{RB-G-matrix-form}) involves $N \times N_e+N_{\text{el}}$ unknowns, where $N \times N_e+N_{\text{el}}< N_e^f+N_{\text{el}}^f$.

We will present details of the two sampling strategies used for constructing the reduced mixed GMsFE basis:  basis-oriented cross-validation and  POD.

\section{Strategies  for constructing   reduced mixed GMsFE basis functions}
\label{ssec:sampling algorithm}

We construct the set of snapshots for reduced basis  by computing the mixed GMsFE basis functions for each $\mu\in\Xi_{op}$, and denote them by
\[
\Sigma^{i}=\{\phi_j^i(\mu_k):1\leq j\leq l_i(\mu_k), \forall \mu_k \in \Xi_{op}\}, ~~1\leq i\leq N_e.
\]
To simplify notation, we use the following single-index notation
\[
\Sigma^{i}=\{\phi_n^i:1\leq n\leq n_{\text{snap}}^i\}, ~~1\leq i\leq N_e,
\]
where $n_{\text{snap}}^i=\sum_{j=1}^{|\Xi_{op}|}l_i(\mu_j),~ \mu_j\in \Xi_{op}$.  Let
$\Sigma=\{\Sigma^{i}\}_{i=1}^{N_e}$
be the set of the snapshots for reduced mixed GMsFE basis and $n_{\text{snap}}=\min_{1\leq i\leq N_e}\{n_{\text{snap}}^i\}$.

\subsection{Basis-oriented cross-validation method for reduced mixed GMsFE basis}
\label{BOCV sampling algorithm}

   The cross-validation method is devoted to selecting the optimal parameters for multiscale basis in \cite{jiang2016reduced}. In the paper, we will  extend
     the idea  to identify  reduced mixed GMsFE basis functions from the snapshots $\Sigma$. Thus we call it as basis-oriented cross-validation (BOCV).

     Let $\Xi_{\text{validate}} \subset \Gamma$ be a given validation set, and $\varepsilon^{*}$ a tolerance for the stopping criterion for the greedy algorithm. The main steps of the basis-oriented cross-validation method are as follows:\\
   $\bullet$ Split the snapshots $\Sigma$ into $n_{\text{snap}}$ disjoint subset $\digamma^n=\{\phi_n^i: 1\leq i\leq N_e\},$ $n=1,...,n_{\text{snap}};$ \\
   $\bullet$ For each subset $\digamma^n$, we solve the equation (\ref{RB-G-matrix-form}) in the corresponding space $\text{span}\{\phi_n^i: 1\leq i\leq N_e\}$ for each $\mu \in \Xi_{\text{validate}}$, and compute the mean errors for the validation set;\\
   $\bullet$ Detect the subset $\digamma^n$ with  the minimum mean error, and update the reduced pace with $V_H^{N,CV}=V_H^{N-1,CV}\cup \text{span}\{\phi: \phi\in \digamma^n\}$. Here we initialize the $V_H^{0,CV}={\O}$;\\
   $\bullet$ Repeat the process until the prescribed error tolerance $\varepsilon^{*}$ is satisfied.

   We describe the details of BOCV in Algorithm \ref{algorithm-BOCV}. Here $M_C$ is the number of local basis functions for velocity.
  \begin{algorithm}
  \caption{ Basis-oriented cross-validation for reduced mixed GMsFEM}
     \textbf{Input}: Some snapshots $\Sigma=\{\digamma^n\}_{n=1}^{n_{\text{snap}}}$, a validating set $\Xi_{\text{validate}} \subset \Gamma$, $n_{\text{snap}}$ and \\a tolerance $\varepsilon^{*}$\\
     \textbf{Output}:  the reduced mixed GMsB space: $ V_H^{N,CV}$  \\
      ~1:~$~$\textbf{Initialization}: $N=1$, $V_H^{0,CV}={\O}$; \\
      ~2:~$~$\textbf{for} $j=1:n_{\text{snap}}$\\
      ~3: ~~~~~$e^{j}(\mu)=v_h(\mu)-v_H^{N,j}(\mu),$ \text{where}~ $v_H^{N,j}(\mu) ~ \text{solves} ~(\ref{RB-discre-mix-eq})~  \forall ~u \in V_H^{N,j} $, \\
      $~~~~~~~~~$-\text{and} $ V_H^{N,j}=V_H^{N-1,CV}\cup \text{span}\{\phi:\phi\in \digamma^j\};$\\
      ~4:$~~~~~\mathcal{E}_{\text{mean}}^N(j)=\text{mean}_{\mu\in\Xi_{\text{validate}}}\|e^{j}(\mu)\|_{V}$; \\
      ~5:~$~$\textbf{end for}\\
      ~6:~~$~\text{index}=\arg\min_j \mathcal{E}_{\text{mean}}^N(j);$ \\
      ~7:~~$~V_H^{N,CV}=V_H^{N-1,CV}\cup V_H^{N,\text{index}}:=\text{span}\{\phi_j^i, ~1\leq j\leq N,~1\leq i\leq N_e\};$ \\
      ~8:~~$~\varepsilon=\max_j \mathcal{E}_{\text{mean}}^N(j);$\\
      ~9:$~~~$update $\Sigma$ with $\Sigma=\Sigma \setminus \digamma^{\text{index}};$\\
      10:~~$n_{\text{snap}}\leftarrow n_{\text{snap}}-1$; \\
      11:~~$N\leftarrow N+1$; \\
      12:~~\textbf{if} $\varepsilon \leq \varepsilon^{*} $\\
      13:$~~~~~~$Go back to step 2;\\
      14:~~\textbf{end if} \\
      15:~~$M_C=N$;
  \label{algorithm-BOCV}
\end{algorithm}

According to Algorithm \ref{algorithm-BOCV}, we get the  reduced mixed GMsB space using BOCV
\[
V_H^{M_C,CV}:=\text{span}\{\varphi_{i}^n: 1\leq i\leq N_e,1\leq n\leq N_{M_C}\}.
\]
To obtain a set of $(\cdot,\cdot)_V$-orthonormal basis functions, we apply the Gram-Schmidt process to the set $\{\varphi_{i}^n: 1\leq i\leq N_e,1\leq n\leq N_{M_C}\}$ in the $(\cdot,\cdot)_V$ inner product, denote them by
 \begin{eqnarray}
\label{MsRB-basis-cv}
\big\{\psi_{i}: 1\leq i\leq M_C \times N_e\big\}.
\end{eqnarray}

\subsection{POD method for reduced mixed GMsFE basis}
\label{POD-MsRB}

POD can be used to construct a low rank approximation for a Hilbert space (ref. \cite{sv08}).
Now we consider to use POD for sampling reduced mixed GMsFE basis functions.  The POD multiscale basis can be constructed based on the given local snapshots  $\Sigma^i$.

For each $1\leq i\leq N_e$, we consider a set of snapshots: $\Sigma^i=\{\phi_n^i:1\leq n\leq n_{\text{snap}}^i\}$, which can be written as
\begin{eqnarray*}
\label{eigen-eq}
\phi_n^i=\sum_{k=1}^{N_e^{f,i}}y_{ki}^n \xi_{k}=\vec{\xi}\vec{y}_{i}^n, ~1\leq n\leq n_{\text{snap}}^i,
\end{eqnarray*}
 where $\vec{\xi}=[\xi_{1},\cdot\cdot\cdot,\xi_{N_e^{f,i}}]$ is the set of mixed FE basis functions for velocity supporting on $w_i$, $N_e^{f,i}$ is the number of the corresponding basis functions, and $\vec{y}_{i}^n=[y_{1i}^n,\cdot\cdot\cdot,y_{N_e^{f,i} i}^n]^T\in\mathbb{R}^{N_e^{f,i}}$. The matrix of the snapshot coefficients is defined by
 \[
 Y_i=(y_{ki}^n)=[\vec{y}_{i}^1, \cdots \vec{y}_{i}^n, \cdots \vec{y}_i^{n_{\text{snap}}^i}]
 \in\mathbb{R}^{N_e^{f,i}\times n_{\text{snap}}^i},\quad 1\leq n\leq n_{\text{snap}}^i,\quad 1\leq k\leq N_e^{f,i}.
 \]
 The POD method for sampling reduced mixed GMsFE basis is described  in Algorithm \ref{algorithm-POD}.

 \begin{algorithm}
\caption{POD method for sampling reduced mixed GMsFE basis}
     \textbf{Input}: Some snapshot $\Sigma=\{\Sigma^n\}_{n=1}^{N_e}$, and  $M_{P}$\\
     \textbf{Output}:  the reduced mixed GMsB space: $V_H^{N,POD}$  \\
      $\diamond$:~~For each $1\leq i\leq N_e$,  based on the snapshots $\Sigma=\{\Sigma^n\}_{n=1}^{N_e}$, the POD basis functions\\
      $~~~~~$can be constructed as follows.\\
      ~~~~~$\triangleright$:~Construct a matrix $\aleph$ using  the inner product of the snapshots.\\
      $~~~~~\aleph=((\phi_m^i,\phi_n^i)_V)=Y_i^{T}M_h^i Y_i\in \mathbb{R}^{n_{\text{snap}}^i\times n_{\text{snap}}^i}.$ where $M_h^i=[(\xi_{k},\xi_{k'})_X]\in\mathbb{R}^{N_e^{f,i}\times N_e^{f,i}}$.\\
      ~~~~~$\triangleright$:~Compute $M_{P}$ eigenvectors of $\aleph$ corresponding the first $M_{P}$ largest eigenvalues, i.e.,\\
      $~~~~~\aleph z_j=z_j\lambda_j, j=1,\cdot\cdot\cdot,M_{P},$ where $z_j\in\mathbb{R}^{n_{\text{snap}}^i}$ and $\lambda_1\geq \lambda_2,\cdot\cdot\cdot,\lambda_{M_{P}}$.\\
      ~~~~~$\triangleright$:~For any $1\leq j\leq M_P$,  the POD basis functions $\psi_{i}^j$ are given by:\\
      $~~~~~\psi_{i}^j(x)=\frac{1}{\sqrt{\lambda_j}}\sum_{n=1}^{n_{\text{snap}}^i}(z_j)_{n}\phi_n^i.$\\
      $\diamond$:~~We get the reduced mixed GMsB space using POD:~$V_H^{M_P,POD}:=\text{span}\{\psi_{i}^j: 1\leq i\leq N_e,1\leq j\leq M_P\}.$
  \label{algorithm-POD}
\end{algorithm}

\subsection{Selecting $\Xi_{op}$ by a greedy algorithm}
\label{greedy selection}

Let $\Xi_{\text{train}}$ be a training set, which is a collection of a finite number of samples in $\Gamma$. Typically the training set is chosen by Monte Carlo methods. It is required that the samples in $\Xi_{\text{train}}$ are sufficiently scattered in the domain $\Gamma$.   We apply a greedy algorithm to get a few optimal samples from the training set, and denote the set of the optimal samples by $\Xi_{op}$.

We use $\Xi_{op}$ instead of $\Xi_{\text{train}}$ to construct the set of snapshots for reduced mixed GMsFEM, which can significantly improve computation efficiency since $|\Xi_{op}|\ll |\Xi_{\text{train}}|$. We consider an error bound for mixed GMsFEM. There are two basic ingredients of the error bound: residual error and stability information of the corresponding bilinear form $a(\cdot,\cdot;\mu)$ and $b(\cdot,\cdot;\mu)$.

 First we consider the residual error for  mixed GMsFEM, which is  important  for posteriori analysis.
 Let
 \[
e_v(\mu): = v_h(\mu)-v_H^N (\mu), \quad  e_p(\mu): = p_h(\mu)-p_H^N (\mu).
\]
 By equation (\ref{mix-FEM-dis}), we get
\begin{eqnarray*}
\label{deduce-equation}
\begin{cases}
\begin{split}
a\big(e_v(\mu),u; \mu\big)-b\big(u,  e_p(\mu); \mu\big)&=-a\big(v_H^N (\mu),u; \mu\big)+b\big(u,p_H^N (\mu); \mu\big), ~~~\forall ~u \in V_h, \\
b\big(e_v(\mu), q; \mu\big)&=l(q)-b\big(v_H^N (\mu), q; \mu\big), ~~~\forall ~q \in Q_h.
\end{split}
\end{cases}
\end{eqnarray*}
Let $r_1(u;\mu)\in V_h^{*}$ (the dual space of $ V_h$) and $r_2(q;\mu)\in Q_h^{*}$  be the residual
\begin{eqnarray*}
\label{Rfa-weak-eq}
\begin{cases}
\begin{split}
r_1(u;\mu):&=-a\big(v_H^N (\mu),u; \mu\big)+b\big(u,p_H^N (\mu); \mu\big), ~~~\forall ~u \in V_h,\\
r_2(q;\mu):&=l(q)-b\big(v_H^N (\mu), q; \mu\big), ~~~\forall ~q \in Q_h.
\end{split}
\end{cases}
\end{eqnarray*}
Then  we get
\begin{eqnarray}
\label{error-residual-eq}
\begin{cases}
\begin{split}
a\big(e_v(\mu),u; \mu)-b\big(u,e_p(\mu); \mu\big)&=r_1(u;\mu), ~~~\forall ~u \in V_h,\\
l(q)-b\big(e_v(\mu), q; \mu\big)&=r_2(q;\mu), ~~~\forall ~q \in Q_h.
\end{split}
\end{cases}
\end{eqnarray}
By Riesz representation theory, there exist $\hat{e}_v(\mu)\in V_h$ and $\hat{e}_p(\mu)\in Q_h$ such that
\begin{eqnarray}
\label{Riesz-weak-eq}
\begin{cases}
\begin{split}
\big(\hat{e}_v(\mu),u\big)_{V}&=r_1(u;\mu),  \quad \forall u\in V_h,\\
\big(\hat{e}_p(\mu),q\big)_{Q}&=r_2(q;\mu),  \quad \forall q\in Q_h.
\end{split}
\end{cases}
\end{eqnarray}
Then we can rewrite the error residual equation (\ref{error-residual-eq}) as
\begin{eqnarray*}
\label{Riesz error residual equation}
\begin{cases}
\begin{split}
a\big(e_v(\mu),u; \mu)-b\big(u,e_p(\mu); \mu\big)&=\big(\hat{e}_v(\mu),u\big)_{V},  \quad \forall u\in V_h,\\
l(q)-b\big(e_v(\mu), q; \mu\big)&=\big(\hat{e}_p(\mu),q\big)_{Q},  \quad \forall q\in Q_h.
\end{split}
\end{cases}
\end{eqnarray*}
Consequently,  the dual norm of the residual  $r_1(u; \mu)$ and $r_2(q; \mu)$  can be evaluated through the Riesz representation,
\begin{eqnarray}
\label{dual-norm-residual}
\begin{cases}
\begin{split}
\|r_1(u;\mu)\|_{V_h^{*}}:=\sup_{u\in V_h}\frac{r_1(u;\mu)}{\|u\|_{V}}=\|\hat{e}_v(\mu)\|_{V},\\
\|r_2(q;\mu)\|_{Q_h^{*}}:=\sup_{q\in Q_h}\frac{r_2(q;\mu)}{\|q\|_{Q}}=\|\hat{e}_p(\mu)\|_{Q}.
\end{split}
\end{cases}
\end{eqnarray}

Secondly, we need a lower bound  and an upper bound for the continuity constant in (\ref{bounded-a}) and coercivity constant in (\ref{coercive-a}),
\[
0<\gamma_{LB}(\mu)\leq\gamma(\mu)\leq\gamma_{UB}(\mu), ~~~ \forall \mu \in \Gamma,
\]
\[
0<\alpha_{LB}(\mu)\leq\alpha(\mu)\leq\alpha_{UB}(\mu), ~~~ \forall \mu \in \Gamma,
\]
and for the inf-sup constant in (\ref{inf-sup condition}),
\[
0<\beta_{LB}(\mu)\leq\beta(\mu)\leq\beta_{UB}(\mu), ~~~ \forall \mu \in \Gamma.
\]
Based on these bounds,  we define an error estimator \cite{ gerner2012certified, rhp08} for the solution of equation (\ref{error-residual-eq}) by
\begin{eqnarray}
\label{erroe estimator-v}
\Delta_{N}^v(\mu):=\frac{\|\hat{e}_v(\mu)\|_{V}}{{\alpha_{LB}(\mu)}}+\bigg(1+\frac{\gamma_{UB}(\mu)}{\alpha_{LB}(\mu)}\bigg)\frac{\|\hat{e}_p(\mu)\|_{Q}}{{\beta_{LB}(\mu)}},
\end{eqnarray}
\[
\Delta_{N}^p(\mu):=\frac{\|\hat{e}_v(\mu)\|_{V}}{{\beta_{LB}(\mu)}}+\frac{\gamma_{UB}(\mu)}{\beta_{LB}(\mu)}\Delta_{N}^v(\mu).
\]

An efficient method for computing $\alpha_{LB}(\mu)$, $\gamma_{UB}(\mu)$ and $\beta_{LB}(\mu)$ is the Successive Constraint Method (ref. \cite{ gerner2012certified, pgsa07, rhp08}).
Note that $\Delta_{N}^v(\mu)$ and $\Delta_{N}^p(\mu)$ are the upper bounds for the errors $\|e_v(\mu)\|_{V}$ and $\|e_p(\mu)\|_{Q}$ such that
\[
\|e_v(\mu)\|_{V}\leq  \Delta_{N}^v(\mu), ~~~\|e_p(\mu)\|_{Q}\leq  \Delta_{N}^p(\mu), ~~~ \forall \mu \in \Gamma.
\]
In this paper, the velocity is what we are particularly interested in. Thus we introduce the associated effectivity with the error estimator $\Delta_{N}^v(\mu)$,
\[
\eta_{N}^v(\mu):=\frac{\Delta_{N}^v(\mu)}{\|e_v(\mu)\|_{V}}.
\]
The effectivity is a measure of the quality of the proposed estimator. Following the references \cite{gerner2012certified, quarteroni2011certified, rhp08}, the effectivity satisfies
\[
1\leq \eta_{N}^v(\mu)\leq \frac{\gamma(\mu)}{\alpha_{LB}(\mu)}+\frac{\gamma_b}{\alpha_{LB}(\mu)}\frac{\|e_p(\mu)\|_{Q}}{\|e_v(\mu)\|_{V}}+\bigg(1+\frac{\gamma_{UB}(\mu)}{\alpha_{LB}(\mu)}\bigg)\frac{\gamma_b}{{\beta_{LB}(\mu)}}, ~~~\forall \mu\in \Gamma.
\]
In the process of iteration, we just update the reduced mixed GMsFEM basis for velocity, the pressure space is constructed by piecewise constant functions on the coarse grid. Thus the inequality  $\|e_v(\mu)\|_{V}\leq\|e_p(\mu)\|_{Q}$ holds in general.
Then the associated effectivity  satisfies
\[
1\leq \eta_{N}^v(\mu)\leq \frac{\gamma(\mu)+\gamma_b}{\alpha_{LB}(\mu)}+\bigg(1+\frac{\gamma_{UB}(\mu)}{\alpha_{LB}(\mu)}\bigg)\frac{\gamma_b}{{\beta_{LB}(\mu)}}, ~~~\forall \mu\in \Gamma.
\]

We note that the $\hat{e}(\mu)$ is related to $r(v;\mu)$ by the equation (\ref{dual-norm-residual}).
By (\ref{GMsRB solution}) and (\ref{affinely-ag}),  the residual can be expressed by
%\begin{eqnarray}
%\label{Qresidual-eq}
%\begin{split}
%r_1(u;\mu)&=-a\big(v_H^N (\mu),u; \mu)+b\big(u,p_H^N (\mu); \mu) \\
%&=-\sum_{i=1}^{N_e \times N}v_i^N(\mu)a(\psi_{i},u; \mu)+\sum_{r=1}^{N_{\text{el}}} p_r^N(\mu)b(u, \eta_{r})\\
%&=-\sum_{i=1}^{N_e \times N}\sum_{j=1}^{m_{a}}k^{j}(\mu)v_i^N(\mu)a^j(\psi_{i},u)+\sum_{r=1}^{N_{\text{el}}} p_r^N(\mu) b(u, \eta_{r}), \quad \forall u\in V_h.\\
%r_2(q;\mu)&=l(q)-\sum_{i=1}^{N_e \times N}v_i^N(\mu)b(\psi_{i}, q) ~~~\forall ~q \in Q_h.
%\end{split}
%\end{eqnarray}
\begin{eqnarray}
\label{Qresidual-eq}
\begin{split}
r_1(u;\mu)&=-a\big(v_H^N (\mu),u; \mu\big)+b\big(u,p_H^N (\mu); \mu\big) \\
&=-\sum_{i=1}^{N_e \times N}\sum_{j=1}^{m_{a}}k^{j}(\mu)v_i^N(\mu)a^j(\psi_{i},u)+\sum_{r=1}^{N_{\text{el}}} p_r^N(\mu) b(u, \eta_{r}), \quad \forall u\in V_h,\\
r_2(q;\mu)&=l(q)-\sum_{i=1}^{N_e \times N}v_i^N(\mu)b(\psi_{i}, q) ~~~\forall ~q \in Q_h.
\end{split}
\end{eqnarray}

By (\ref{Qresidual-eq}) and (\ref{Riesz-weak-eq}), we have
\begin{eqnarray*}
\label{Riesz Qerror estimators}
\begin{cases}
\begin{split}
\big(\hat{e}_v(\mu),u\big)_{V}&=-\sum_{i=1}^{N_e \times N}\sum_{j=1}^{m_{a}}k^{j}(\mu)v_i^N(\mu)a^j(\psi_{i},u)+\sum_{r=1}^{N_{\text{el}}}p_r^N(\mu) b(u, \eta_{r}), \quad \forall u\in V_h,\\
\big(\hat{e}_p(\mu),q\big)_{Q}&=l(q)-\sum_{i=1}^{N_e \times N}v_i^N(\mu)b(\psi_{i}, q) ~~~\forall ~q \in Q_h.
\end{split}
\end{cases}
\end{eqnarray*}
This implies that
\begin{eqnarray}
\label{Riesz-error-estimators1}
\hat{e}_v(\mu)=\sum_{i=1}^{N_e \times N}\sum_{j=1}^{m_{a}}k^{j}(\mu)v_i^N(\mu)\mathcal{L}_{i}^{p}+\sum_{r=1}^{N_{\text{el}}}p_r^N(\mu)\mathcal{X}^r,
\end{eqnarray}
\begin{eqnarray}
\label{Riesz-error-estimators2}
\hat{e}_p(\mu)=\mathcal{C}+\sum_{i=1}^{N_e \times N}v_i^N(\mu)\mathcal{X}_i,
\end{eqnarray}
where $\mathcal{L}_{i}^{j}$  is the Riesz representation of $a^{j}(\psi_{i},u)$, i.e., $(\mathcal{L}_{i}^{j},u)_{V}=-a^{j}(\psi_{i},u)$ for any $u\in V_h$, $\mathcal{X}^{r}$  is the Riesz representation of $b(u, \eta_{r})$, i.e., $(\mathcal{X}^{r},u)_{V}=b(u, \eta_{r})$ for any $u\in V_h$, $\mathcal{X}_{i}$  is the Riesz representation of $b(\psi_{i}, q)$, i.e., $(\mathcal{X}_{r},q)_{Q}=-b(\psi_{i}, q)$ for any $q\in Q_h$, and $\mathcal{C}$ is the Riesz representation of $l$, i.e., $(\mathcal{C},q)_{X}=l(q)$ for any $q\in Q_h$, where $1\leq j\leq m_{a}$, $1\leq i\leq N\times N_e$, and $1\leq r\leq N_{\text{el}}$.
The equations (\ref{Riesz-error-estimators1}) and (\ref{Riesz-error-estimators2}) give rise to
\begin{eqnarray}
\label{error-norm-est}
\begin{split}
\|\hat{e}_v(\mu)\|_{V}^2=&\sum_{r=1}^{N_{\text{el}}}\sum_{r'=1}^{N_{\text{el}}}p_r^N(\mu)p_{r'}^N(\mu)(\mathcal{X}^r,\mathcal{X}^{r'})_{V}+\sum_{i=1}^{N_e \times N}\sum_{j=1}^{m_{a}}k^{j}(\mu)v_i^N(\mu)\\
&\times\bigg\{2\sum_{r=1}^{N_{\text{el}}}p_r^N(\mu)(\mathcal{X}^r,\mathcal{L}_{i}^{j})_{V}+\sum_{i'=1}^{N_e \times N}\sum_{j'=1}^{m_{a}}k^{j'}(\mu)v_{i'}^N(\mu)
(\mathcal{L}_{i}^{j},\mathcal{L}_{i'}^{j'})_{V}\bigg\},\\
\|\hat{e}_p(\mu)\|_{Q}^2=&(\mathcal{C},\mathcal{C})_{Q}+2\sum_{i=1}^{N\times N_e}v_i^N(\mu)(\mathcal{C},\mathcal{X}_i)_{Q}+\sum_{i=1}^{N\times N_e}\sum_{i'=1}^{N\times N_e}v_{i'}^N(\mu)v_i^N(\mu)(\mathcal{X}_i,\mathcal{X}_{i'})_{Q}.
\end{split}
\end{eqnarray}

To efficiently compute $\|\hat{e}_v(\mu)\|_{V}$ and $\|\hat{e}_p(\mu)\|_{Q}$, we apply an offline-online procedure. In the offline stage we compute and store the parameter-independent quantities. In particular, we compute $\mathcal{C}$, $\mathcal{X}^{r}$, $\mathcal{X}_{i}$ and $\mathcal{L}_{i}^{j}$, and store $(\mathcal{X}^r,\mathcal{X}^{r'})_{V}$, $(\mathcal{X}^r,\mathcal{L}_{i}^{j})_{V}$, $(\mathcal{L}_{i}^{j},\mathcal{L}_{i'}^{j'})_{V}$, $(\mathcal{C},\mathcal{C})_{Q}$, $(\mathcal{C},\mathcal{X}_i)_{Q}$ and $(\mathcal{X}_i,\mathcal{X}_{i'})_{Q}$, where $1\leq i,i'\leq N\times N_c$, $1\leq j,j'\leq m_{a}$,  $1\leq r,r'\leq N_{\text{el}}$.
In the online stage, for any $\mu$, we compute $v_i^N(\mu)$ ($1\leq i\leq N\times N_e$) and use (\ref{error-norm-est}) to compute $\|\hat{e}_v(\mu)\|_{V}$ and $\|\hat{e}_p(\mu)\|_{Q}$.

In  summary,  we describe the greedy algorithm for selecting $\Xi_{op}$ in Algorithm \ref{algorithm-Greedy}.
\begin{algorithm}[hbtp]
\caption{ Greedy algorithm for selecting  $\Xi_{op}$}
     \textbf{Input}: A training set $\Xi_{\text{train}} \subset \Gamma$, $\mu_1\in \Xi_{\text{train}}$, and $N_p$ \\
     \textbf{Output}:  The optimal set $\Xi_{op}\in \Xi_{\text{train}}$\\
      ~1:~~Initialize $N=1$, $\Xi_{op}=\{\mu_1\}$ and $\mu_{N}=\mu_{1}$;\\
      ~2:~~Construct the snapshots $\Sigma$ for reduced mixed GMsFEM by computing the mixed \\
      $~~~~~$ GMsFE basis for each $\mu\in\Xi_{op}$, and construct reduced basis space $V_H^{N_p}$ by POD;\\
      ~3:~~Update $\Xi_{\text{train}}$ with $\Xi_{\text{train}}=\Xi_{\text{train}}\setminus \mu_{N};$\\
      ~4:~~For each $\mu\in \Xi_{\text{train}}$, compute the reduced basis approximation $v_H^{N_p}$ by (\ref{RB-G-matrix-form}) in $V_H^{N_p}$.\\
      ~5:~~For each $\mu\in \Xi_{\text{train}}$, evaluate the error estimator $\Delta_{N_p}^v(\mu)$ by (\ref{erroe estimator-v}).\\
      ~6:~~Choose $\mu_{N+1}=\arg\max_{\mu\in\Xi_{\text{train}}}\Delta_{N_p}^v(\mu)$, and set $\varepsilon_N=\max_{\mu\in\Xi_{\text{train}}}\Delta_{N_p}^v(\mu)$;\\
      ~7:~~Update $\Xi_{op}$ with $\Xi_{op}=\Xi_{op}\bigcup\{\mu_{N}\}$;\\
      ~8:~~$N\leftarrow N+1$\\
      ~9:~~Return to Step 2 if $\varepsilon_N \leq \varepsilon_{N-1}$, otherwise \textbf{terminate}.
\label{algorithm-Greedy}
\end{algorithm}

\begin{rem}
 We can use the  cross-validation method in \cite{jiang2016reduced} to choose the first optimal sample  $\mu_1$ from the training set $\Xi_{\text{train}}$.
 This choice can improve the accuracy.
\end{rem}
%%%%%%%%%%%%%%%%%%%%%%%%%%%%%%%%%%%%%%%%%%%%%%%%%%%%%%%%%%
 %Our numerical results show that using the optimal $\mu_1$ can give better approximation than randomly chosen $\mu_1$.

\section{Reduced model representation method}
\label{ssec:Online-online}

The reduced order model defined in (\ref{RB-G-matrix-form}) for mixed GMsFEM is a linear algebraic system with $N \times N_e+N_{\text{el}}$ unknowns. Although  the reduced
model needs much less computation effort than original full order model, it may be  not a very small-scale problem because  the reduced order model defined in (\ref{RB-G-matrix-form}) for mixed GMsFEM involves  $N \times N_e+N_{\text{el}}$ unknowns.
In order to significantly improve the online computation efficiency,
we want to get a representation like the form  (\ref{separating presentation}) for the solution and model outputs.
We call this by reduced model representation method, which allows online computation by direct evaluation.
 We employ  LSMOS and STAOMP presented in Section \ref{ssec:PGDDF-LSMOS-STAOMP}   to get the reduced model representation. Let $\Xi_{t}$ be a collection of a finite number of samples in $\Gamma$, and $\{v(x,\mu_i), p(x,\mu_i)\}_{i=1}^{n_t}$ are the snapshots solved by the reduced model (\ref{RB-G-matrix-form}), where $n_t=|\Xi_{t}|$.

 Algorithm \ref{algorithm-LSMOS} outlines the reduced model representation method by LSMOS.
\begin{algorithm}[hbtp]
\caption{ Least-squares method of snapshots for reduced model representation}
     \textbf{Input}: A training set $\Xi_{t} \subset \Gamma$ and a tolerance $\varepsilon^\text{on}$\\
     \textbf{Output}: Reduced model representation $G(x,\mu)\approx \sum_{i=1}^{M} c_i \zeta_i(\mu)v_i(x)$\\
      ~1:~~Compute the snapshots $\{v(x,\mu_i)\}_{i=1}^{n_{n_t}}$ by (\ref{RB-G-matrix-form}) for all $\mu_i\in \Xi_{t}$ and construct \\
      $~~~~~$the covariance
       matrix \textbf{C};\\
      ~2:~~Solve the eigenvalue problem and determine $M$ such that $\frac{\sum _{i=1}^{M}\hat{\lambda}_i}{\sum _{i=1}^{n_t}\hat{\lambda}_i}\leq\varepsilon^\text{on};$\\
      ~3:~~Assemble $\textbf{A}$ based on GPC basis functions and $\Xi_{t}$ by (\ref{least square A});\\
      ~4:~~Construct the functions $\{g_i(x)\}_{i=1}^{M}$ by (\ref{KLE basis}), for each $i=1,..,M$, assemble \textbf{F}\\
      $~~~~~$ by (\ref{least square F});\\
      ~5:~~For each $i=1,..,M$, solve problem (\ref{least square problem}) by least square procedure to obtain $\textbf{d}$,\\
      $~~~~~$ and then get $\zeta_i(\mu)\approx\sum_{i=1}^{M_{g}}d_ip_i(\mu)$;\\
      ~6:~~Return the representation  $G(x,\mu)\approx \sum_{i=1}^{M} \sqrt{\hat{\lambda}_i} \zeta_i(\mu)g_i(x)\approx\sum_{i=1}^{M}\sum_{j=1}^{M_{g}} \sqrt{\hat{\lambda}_i}d_jp_j(\mu) g_i(x)$.
   \label{algorithm-LSMOS}
\end{algorithm}

Algorithm \ref{algorithm-STAOMP} combines reduced mixed GMsFE basis methods and STAOMP together and provides a sparse representation for the model.
STAOMP needs much fewer samples than LSMOS. Thus,
the cardinality of the training set in Algorithm \ref{algorithm-STAOMP} is much less than the cardinality of the training set in Algorithm  \ref{algorithm-LSMOS}.
The whole process for STAOMP based on reduced mixed GMsFE basis methods  is illustrated in Figure \ref{Schema}.
%%%%%%%%%%%%%%%%%%%%%%%%%%%%%%%%%%%%%%%%%%%%
\begin{algorithm}[hbtp]
\caption{ Sparse tensor approximation based on Orthogonal-Matching-Pursuit for reduced model representation}
     \textbf{Input}: A training set $\Xi_{t} \subset \Gamma$, the sample data $\{(x_i,\mu_j)\}\subseteq D \times \Gamma$ ($1\leq i\leq n_x$, $1\leq j\leq n_{\mu}$),
     the number of the optimal basis functions $M$ and a tolerance $\varepsilon^\text{on}$\\
     \textbf{Output}:  Reduced model representation $G(x,\mu)\approx \sum_{i=1}^{\text{Mt}} \textbf{c}(\emph{I}(i)) \Psi_{\emph{I}(i)}(\mu, x);$\\
      ~1:~~Compute the snapshots $\{v(x,\mu_i)\}_{i=1}^{n_t}$ by (\ref{RB-G-matrix-form}) for all $\mu_i\in \Xi_{t}$ and construct the\\
       $~~~~~$optimal basis functions $\{v_j(x)\}_{j=1}^{N}$ for  $\mathcal{H}_N\subseteq \mathcal{H}$;\\
      ~2:~~Choose orthogonal polynomials basis functions $\{p_i(\mu)\}_{i=1}^M$ for $L_2(\Gamma)$, and construct the
      $~~~~~$finite dimensional  approximation space $S_{N\times M}=\{\sum_{i\in I}w_i\Psi_i(x,\mu); ~w_i\in \mathbb{R}\}$\\
      ~3:~~Assemble matrix $\Pi$ and vector $\textbf{b}$ by (\ref{Pi-b-matrix});\\
      ~4:~~Solve the optimization problem (\ref{P0-PROBLEM-discrete}) by Algorithm \ref{algorithm-OMP} with $\varepsilon=\varepsilon^\text{on}$ and get the sparse \\
      $~~~~~$solution $\textbf{c}$ and the solution support $\emph{I}$;\\
      ~5:~~Return the representation $G(x,\mu)\approx \sum_{i=1}^{\text{Mt}} \textbf{c}(\emph{I}(i)) \Psi_{\emph{I}(i)}(\mu, x).$
    \label{algorithm-STAOMP}
\end{algorithm}

\section{Numerical results}
\label{ssec:Numerical result}
In this section, we present a few  examples to illustrate the performance of the proposed reduced mixed GMsFE basis methods
 and make comparisons for different strategies for the model reduction. In Section \ref{Numerical-result1}, we consider an example to illustrate performance of the different reduced mixed GMsFE basis methods for elliptic PDEs with one-dimensional parameters. In Section \ref{Numerical-result2}, we study reduced model representation method for multiscale
 elliptic PDEs with multivariate parameters.
 In Section \ref{Numerical-result3}, we consider a two-phase flow problem in random porous media.  We apply reduced  mixed GMsFE basis methods to the flow equation and integrate STAOMP into the model outputs (water saturation and water-cut) for uncertainty quantification.

\begin{figure}[htbp]
\centering
\subfigure[$\kappa_1$ ]{
    \label{fig:subfig:a}
  \includegraphics[width=3in, height=2.4in]{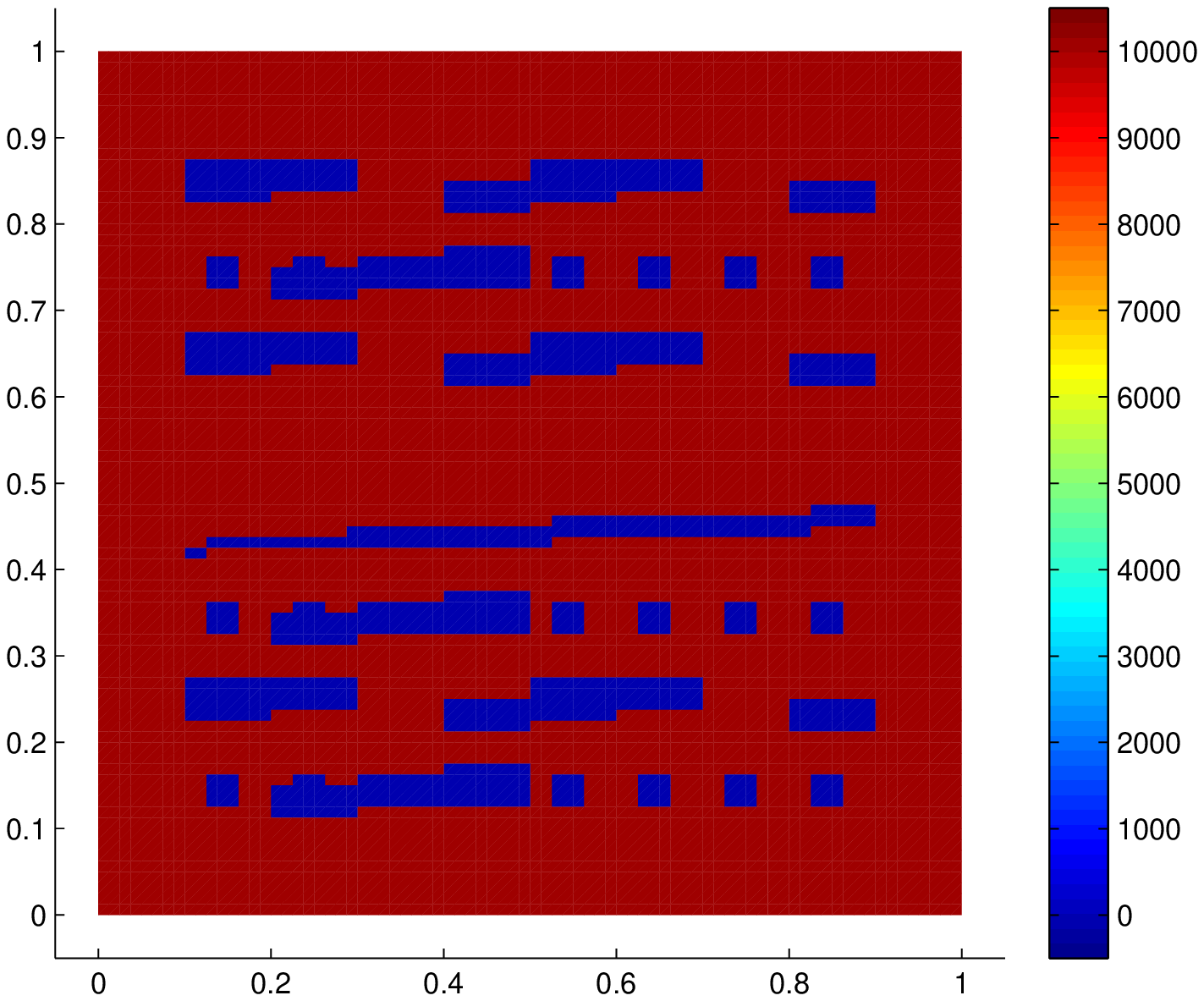}}
%\subfigure[$\kappa_2$]{
%    \label{fig:subfig:b}
%  \includegraphics[width=3in, height=2.4in]{parafield/khc}}
  \subfigure[$\kappa_2$]{
    \label{fig:subfig:b}
  \includegraphics[width=3in, height=2.4in]{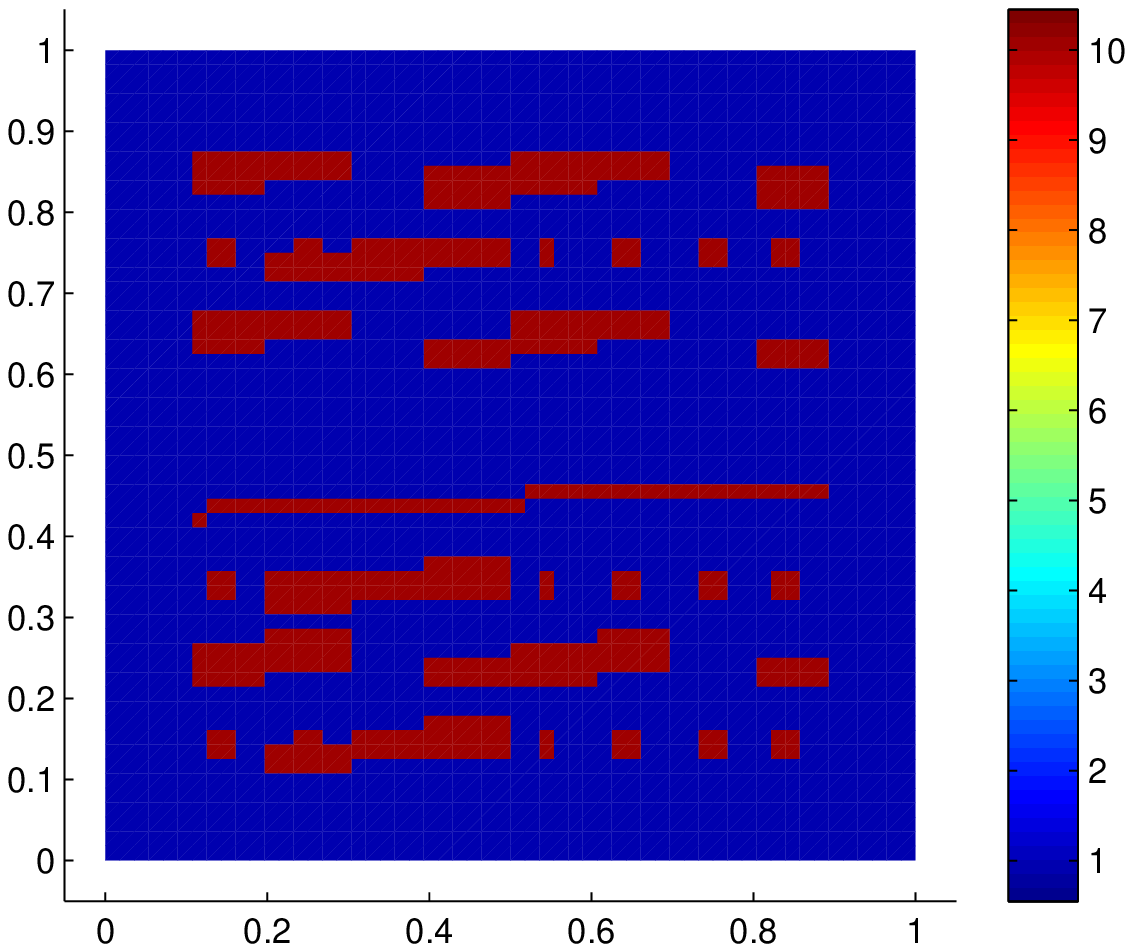}}
  \caption{High-contrast  fields in the numerical experiments.}
  \label{fig-exam}
\end{figure}

 For the numerical examples,  the models are defined in the spatial domain $D=(0,1)^2$. Let $v(x,\mu_i)$ and  $p(x,\mu_i)$ be the reference solutions for velocity and pressure,
 respectively, which are solved by mixed FEM on a fine grid. Let $v_{H}(x,\mu_i)$ and $ p_{H}(x,\mu_i)$ be solved by the reduced mixed GMsFE basis methods.
 Then the relative mean errors for velocity and pressure are defined, respectively, by
\begin{eqnarray}
\label{errors-velocity}
\varepsilon_v=\frac{1}{N}\sum_{i=1}^N\frac{\|v(x,\mu_i)-v_{H}(x,\mu_i)\|_{L^2(D)}}{\|v(x,\mu_i)\|_{L^2(D)}},
\end{eqnarray}
\begin{eqnarray}
\label{errors-pressure}
\varepsilon_p=\frac{1}{N}\sum_{i=1}^N\frac{\|p(x,\mu_i)-p_{H}(x,\mu_i)\|_{L^2(D)}}{\|p(x,\mu_i)\|_{L^2(D)}}.
\end{eqnarray}
 In all of the examples, we will consider the models with  high-contrast  fields, which are depicted in Figure \ref{fig-exam}.

\subsection{Numerical results for reduced mixed GMsFE basis methods}
\label{Numerical-result1}

In this subsection, we consider the elliptic equation in mixed formulation (\ref{mix-eq}),  where  the coefficient function $k(x, \mu)$ is defined by
\[
k(x,\mu)=\frac{10000}{10\sin(20\mu+x_1x_2)+\big(\cos(\mu)+1.2\big)\kappa_1+25},
\]
where $\kappa_1$ is a  high-contrast field and depicted  in Figure \ref{fig-exam} (left), and  the source function $f(x, \mu)$ is
\[
f(x, \mu)=(x_2-0.5)\cos\big(\pi(x_1-0.5)\big).
\]
Here $x:=(x_1, x_2)\in D$ and the random parameter $\mu\sim U(-1,1)$. The diffusion coefficient  $k(x, \mu)$ is oscillating with respect to the random parameter $\mu$ and  highly heterogeneous  with respect to the spatial variable $x$. We apply EIM for $k^{-1}(x, \mu)$  to achieve  offline-online computation decomposition. For the discretization of the spatial domain, we use $80\times80$ uniform fine grid, where the reference solution is computed by the mixed FEM with the lowest Raviart-Thomas FE space. The  mixed GMsFEM is implemented  on the $8\times8$ coarse grid. We choose $|\Xi_{\text{train}}|=200$ parameter values, and select $|\Xi_{op}|=10$ optimal samples from the training set $\Xi_{\text{train}}$ by the greedy algorithm (Algorithm \ref{algorithm-Greedy}). The mixed GMsFEM is used to compute the snapshots, we take $l_i(\mu)=5$ for each coarse block $\omega_i$ and each $\mu\in \Xi_{op}$, thus the number of the snapshot is $n_{\text{snap}}^i=\sum_{j=1}^{n_t}l_i(\mu_j)=50$ for each coarse block $\omega_i$.
To show the efficacy of greedy algorithm for selecting the optimal sample set $\Xi_{op}$, we use two methods to select $\Xi_{op}$: the greedy algorithm and random selection, and then
compare the results.
For the numerical example, we consider the four different reduced GMsFE basis methods: reduced mixed GMsFE  basis method using greedy algorithm and BOCV (RmGMsB-GBOCV), reduced mixed GMsFE  basis method using greedy algorithm and POD (RmGMsB-GPOD), reduced mixed GMsFE  basis method using random selection and BOCV (RmGMsB-RBOCV),  and reduced mixed GMsFE  basis method using random selection and POD (RmGMsB-RPOD). The first two methods use the greedy algorithm to select $\Xi_{op}$, while the last two methods are performed by randomly choosing $10$  parameter samples from $\Xi_{\text{train}}$.

\begin{figure}[htbp]
\centering
  \includegraphics[width=5.5in, height=2.5in]{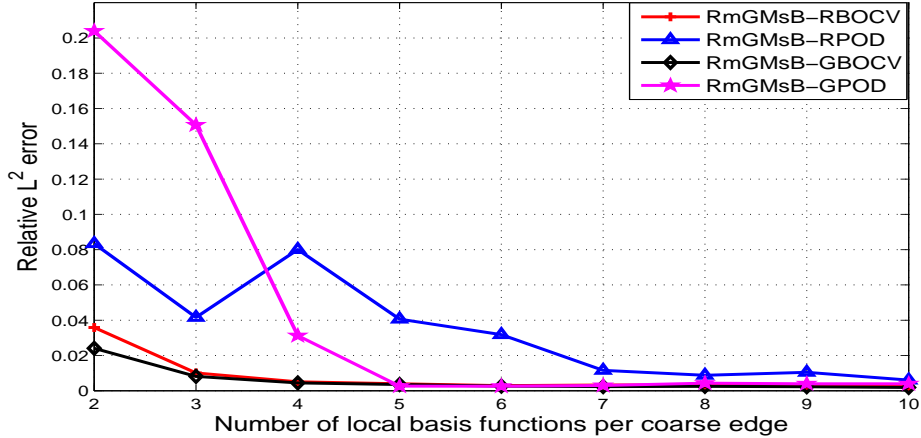}
  \caption{The average relative error for velocity versus number of local basis functions for  RmGMsB-GBOCV, RmGMsB-GPOD, RmGMsB-RBOCV and RmGMsB-RPOD,  $80\times 80$ fine grid, and $8\times 8$ coarse grid. }
  \label{fig1-exam1}
\end{figure}

\begin{figure}[htbp]
\centering
  \includegraphics[width=5in, height=3in]{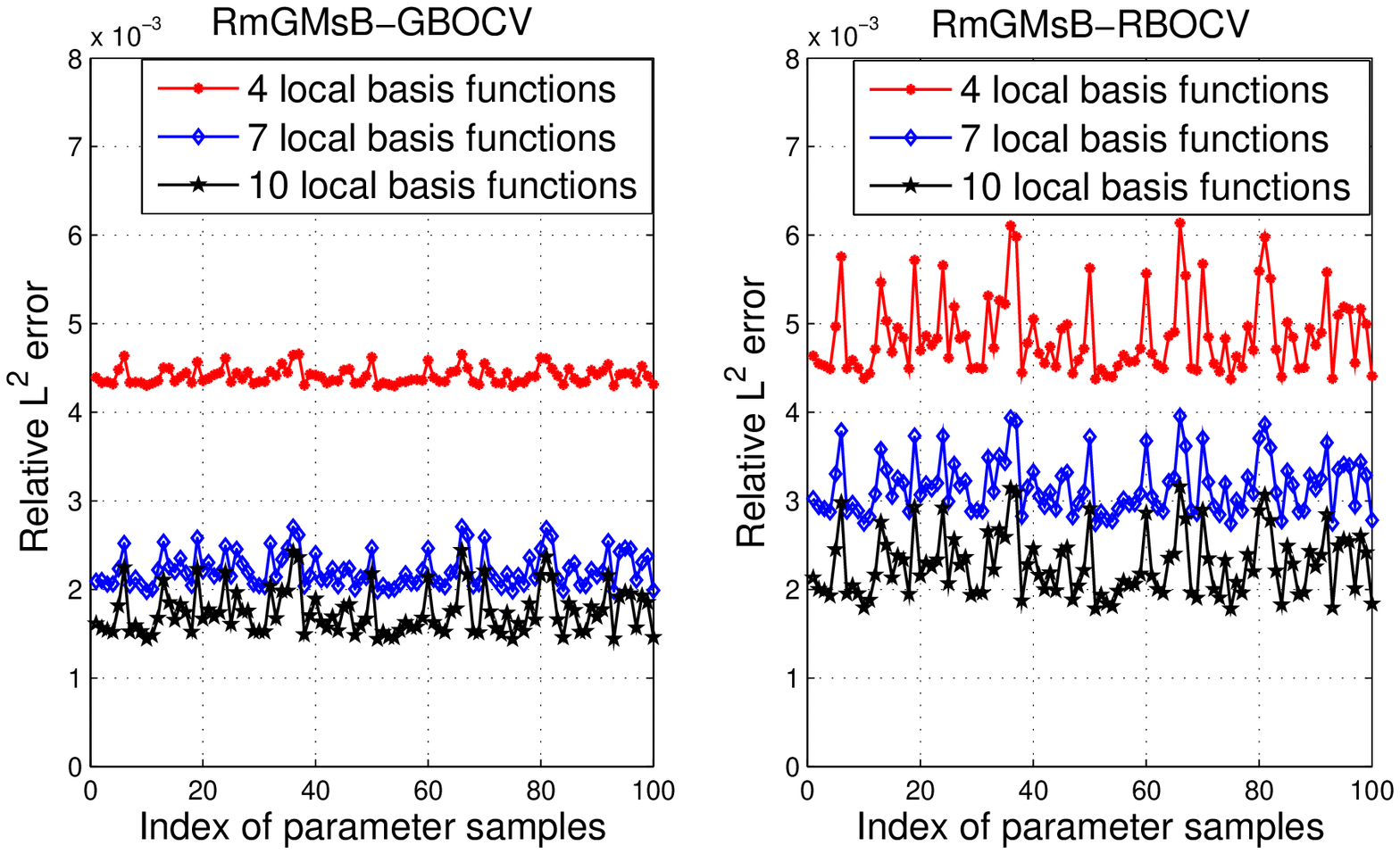}\\
  \includegraphics[width=5in, height=3in]{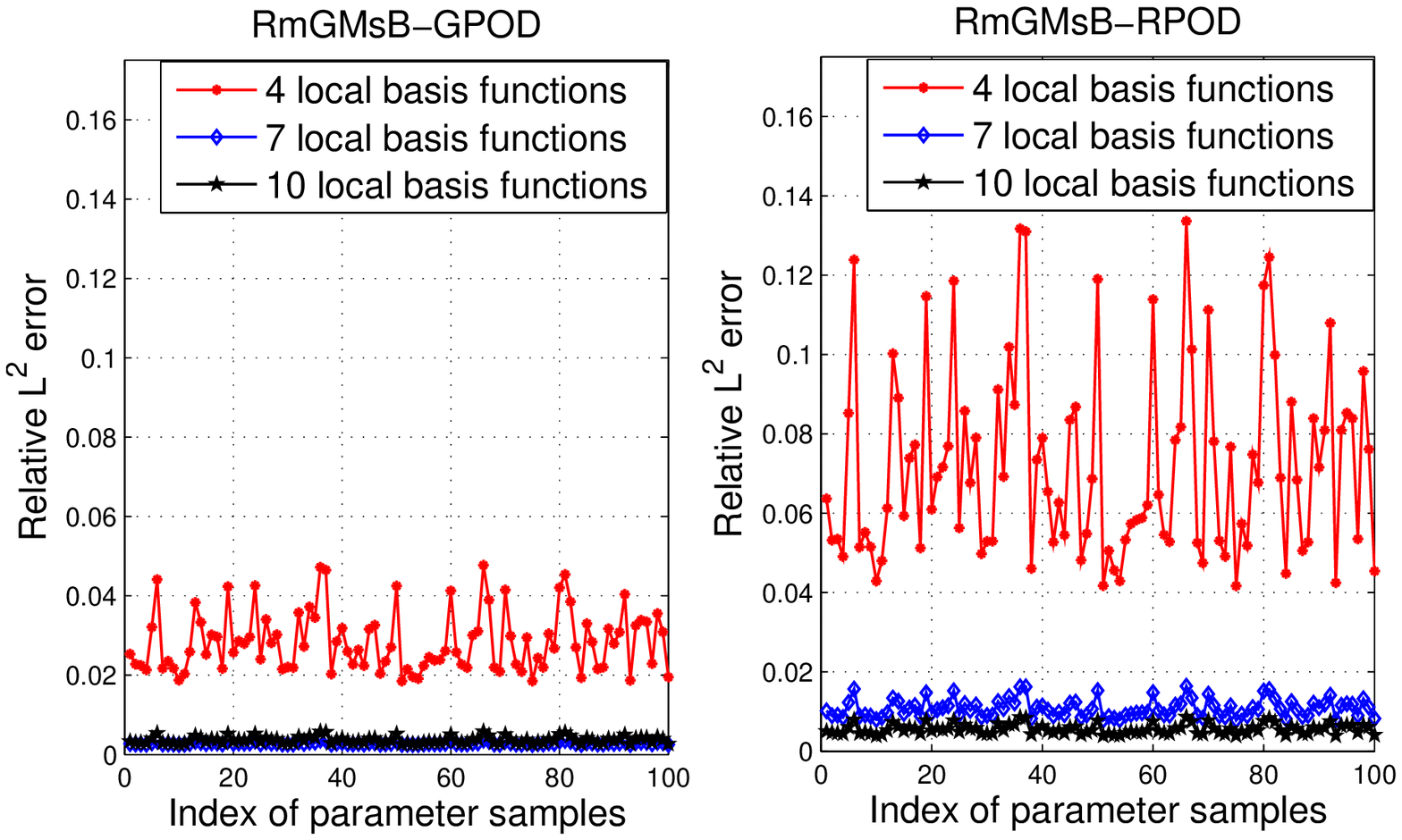}
  \caption{The relative error for $100$ samples by  RmGMsB-GBOCV, RmGMsB-GPOD, RmGMsB-RBOCV and RmGMsB-RPOD,  $80\times 80$ fine grid, and $8\times 8$ coarse grid.}
  \label{fig2-exam1}
\end{figure}

\begin{figure}[htbp]
\centering
  \includegraphics[width=6in, height=3.5in]{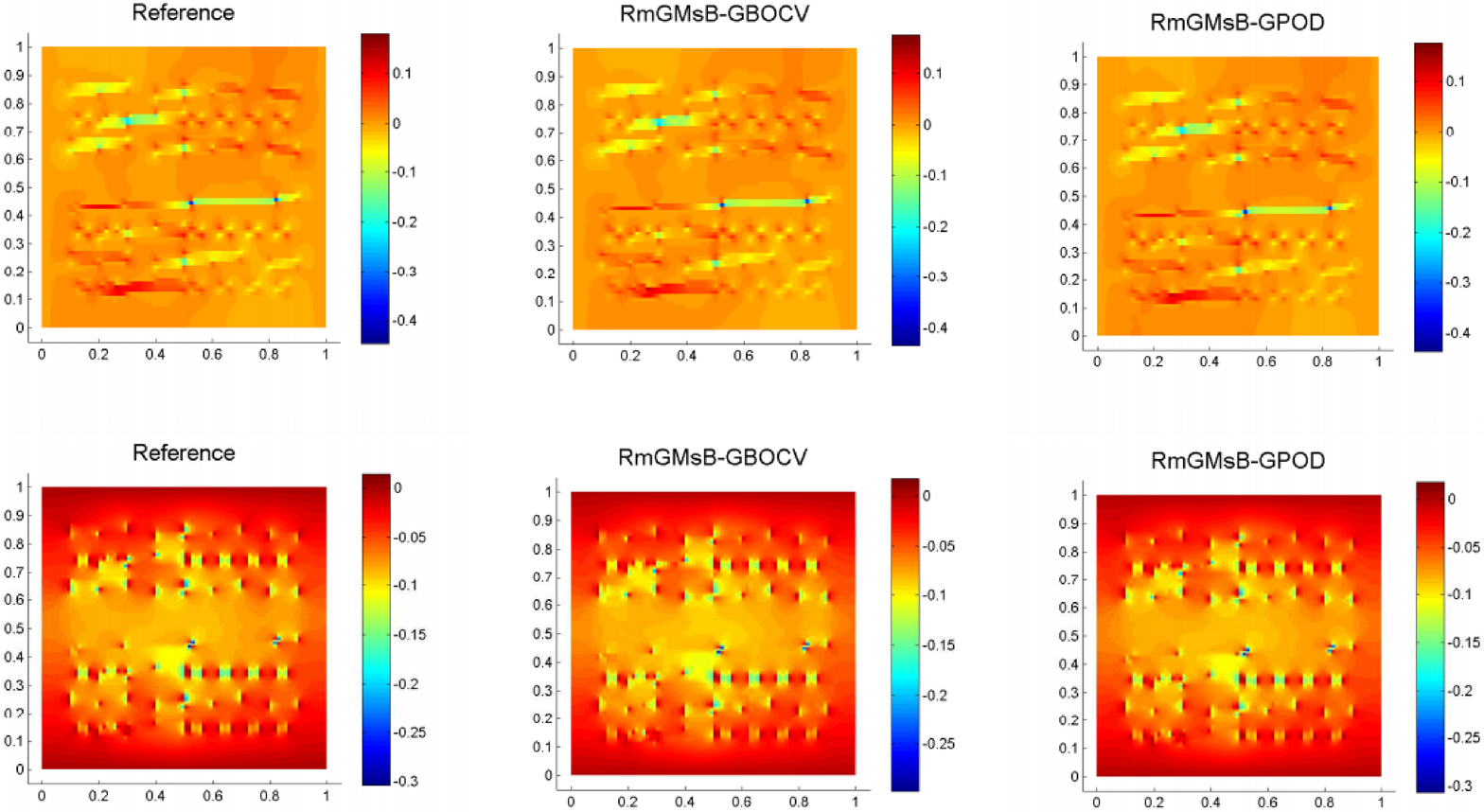}
  \caption{The mean of  velocity profile,  $5$  local GMsFE basis functions on each coarse edge,  the first row are the velocity in x-axis direction and the second row are the
  velocity in y-axis direction.}
  \label{fig4-exam1}
\end{figure}

To compare the approximation accuracy of the four reduced mixed GMsFE  basis methods, we randomly choose $1000$ samples from the parameter space and compute the average of relative  error defined in (\ref{errors-velocity}) and (\ref{errors-pressure}) for the four methods.
In Figure \ref{fig1-exam1}, we depict average relative velocity error versus number of local basis functions for the four reduced mixed GMsFE  basis methods. By the figure we have three observations: (1) as the number of local basis functions increases, the approximation becomes more accurate for all methods; (2) RmGMsB-GBOCV method always achieves better approximation than RmGMsB-RBOCV,  and RmGMsB-GPOD  gives better approximation than RmGMsB-RPOD when the number of local GMsFE basis functions is more than $3$; (3) BOCV approach renders better accuracy than POD approach.
 To visualize the individual errors of the first $100$ samples, we plot the relative errors for the four methods in Figure \ref{fig2-exam1}, which shows that: (1) the error of RmGMsB-GBOCV is much less sensitive to the parameter samples compared with RmGMsB-RBOCV, and the POD approaches have the same situation; (2) BOCV approaches generally gives better approximation than POD approaches. Figure \ref{fig4-exam1} shows the velocity solution profile of the example for the three methods: mixed FEM on fine grid (reference solution), RmGMsB-GBOCV and RmGMsB-GPOD. The figure shows all the velocity  profiles have very good agreement and the reduced mixed GMsFE basis methods provide accurate approximation to the original fine scale model.

%\begin{figure}[htbp]
%\centering
%  \includegraphics[width=2.1in, height=1.8in]{referencev1}
%  \includegraphics[width=2.1in, height=1.8in]{gcvv1}
%  \includegraphics[width=2.1in, height=1.8in]{gpodv1}\\
%  \includegraphics[width=2.1in, height=1.8in]{referencev2}
%  \includegraphics[width=2.1in, height=1.8in]{gcvv2}
%  \includegraphics[width=2.1in, height=1.8in]{gpodv2}
%  \caption{The mean of  velocity profile,  $5$  local GMsFE basis functions on each coarse edge,  the first row are the velocity in x-axis direction and the second row are the
%velocity in y-axis direction.}
%  \label{fig4-exam1}
%\end{figure}

\subsection{Numerical results for reduced model's representation method}
\label{Numerical-result2}
In this section, we consider two examples to illustrate the performance of reduced model representation using LSMOS and STAOMP. In the first example,
we use a multivariate function to highlight the differences between LSMOS and STAOMP.  For the second example, we consider an elliptic PDE with multi-dimensional parameters and apply LSMOS and STAOMP to obtain a reduced model representation for solution.

\subsubsection{Numerical example I: a multivariate function}

Let us consider the function
\[
u(x,\mu):=\sum_{i=1}^2x_i\mu_i+\sin\big(\frac{\pi}{4} (x_1+\frac{1}{3}\prod_{i=1}^3 \mu_i)\big)+\cos\big(\frac{\pi}{4} (x_2+\frac{1}{3}\prod_{i=4}^6 \mu_i)\big),
\]

For the discretization of the spatial domain, we use $50\times50$ uniform grid and will consider the two variable-separation methods: LSMOS and STAOMP.  With regard to LSMOS, we take $\varepsilon^\text{on}=1.3\times 10^{-3}$ in Algorithm \ref{algorithm-LSMOS} and take $900$ parameter samples for snapshots, the representation (\ref{KLE Snapshots}) can be obtained by Algorithm \ref{algorithm-LSMOS} with $M=6$. We use Legendre polynomials with degree up to $N_g=5$ to approximate  $\{\zeta_i(\mu)\}_{i=1}^{M}$ based on least square methods. Hence, the number of the total orthogonal polynomial basis  is $M_g=462$.
 For STAOMP, we take $160$ parameter samples for snapshots, i.e., $n_{t}=160$, and construct $10$ optimal basis functions $\{v_j(x)\}_{j=1}^{6}$ for spatial space, and set $\varepsilon^\text{on}=1.3\times 10^{-3}$ in Algorithm \ref{algorithm-STAOMP}, $160$ parameter samples and $100$ spatial coordinates are selected for snapshots, then we apply  Algorithm \ref{algorithm-STAOMP} to construct the representation $u(x,\mu)\approx \sum_{i=1}^{\text{Mt}} \textbf{c}(\emph{I}(i)) \Psi_{\emph{I}(i)}(\mu, x)$ with $\text{Mt}=41$. Based on these representations, we choose $1000$ samples and compute the average relative error, which is defined as follows,
\begin{eqnarray*}
\label{errors-u}
\varepsilon_u=\frac{1}{N}\sum_{i=1}^N\frac{\|u(x,\mu_i)-\tilde{u}(x,\mu_i)\|_{L^2}}{\|u(x,\mu_i)\|_{L^2}},
\end{eqnarray*}
where $N=1000$ and $\tilde{u}(x,\mu_i)$ is the solution by LSMOS or STAOMP. We list the average relative errors in Table \ref{tab1-exam2} along with the average online CPU time. From the table, we find that: (1) for LSMOS, the approximation error decays fast when  the number of terms  increases; (2) using much fewer terms,  STAOMP can achieve the same accuracy as LSMOS; (3) the magnitude of CPU time for STAOMP is much smaller than that of LSMOS.

\begin{table}[hbtp]
\small
\centering
\caption{ Comparison of relative mean errors $\varepsilon$ and the average online CPU time for different approaches. $\text{Mt}$ is the number of terms retained in  reduced model representation.}
%\vspace*{2pt}
\begin{tabular}{c|c|c|c}
%\hline
\Xhline{1pt}
  Strategies & $\text{Mt}$ & $\varepsilon_u$ & CPU time\\
  \hline
  \multirow{5}{*}{LSMOS}
                & $1\times 462$ & $1.10\times10^{-1}$ & $4.12\times10^{-4}s$ \\
                & $2\times 462$ & $1.33\times10^{-2}$ & $7.80\times10^{-4}s$ \\
                & $3\times 462$ & $1.87\times10^{-3}$ & $1.50\times10^{-3}s$ \\
                & $4\times 462$ & $1.80\times10^{-3}$ & $2.10\times10^{-3}s$ \\
                & $5\times 462$ & $1.70\times10^{-3}$ & $2.30\times10^{-3}s$ \\
 \hline
 \multirow{1}{*}{STAOMP}   & 41 & $1.20\times10^{-3}$ & $9.86\times10^{-5}s$ \\
\Xhline{1pt}
\end{tabular}
\label{tab1-exam2}
\end{table}

%%%%%%%%%%%%%%%%%%%%%%%%%%%%%%%%%%%%%%%
\subsubsection{Numerical example II: an elliptic PDE with multi-dimensional parameters }

In this subsection, we consider the elliptic equation (\ref{mix-eq}) with the source term
\[
f(x)=(x_1+1)\cos(\pi x_2), ~ ~x\in (0, 1)^2.
\]
Let $a(x,\mu)$ be a random field, which is characterized by  a two point exponential  covariance function $\text{cov}[a]$, i.e.,
\begin{eqnarray}
\label{covariance function}
\text{cov}[a](x_1, y_1; x_2, y_2)= \sigma^2 \exp\bigg(-\frac{|x_1-x_2|^2}{2l_x^2}-\frac{|y_1-y_2|^2}{2l_y^2}\bigg),
\end{eqnarray}
where $(x_i, y_i)$ ($i=1,2$) is the spatial coordinate in $D$. Here the variance $\sigma^2=1$, correlation length $l_x=l_y=0.2$. The random coefficient  $a(x, \mu)$ is obtained by truncating a Karhunen-Lo\`{e}ve expansion, i.e.,
\begin{eqnarray}
\label{KLE-TRUNC1}
a(x,\mu):= E[a]+ \sum_{i=1}^{N} \sqrt{\gamma_i}b_i(x) \mu_i,
\end{eqnarray}
where $E[a]=6$ and the random vector $ \mu:=(\mu_1, \mu_2, ...,\mu_{N})\in \mathbb{R}^{N}$. Each $\mu_i$  ($i=1,\cdots, N$) is uniformly distributed in the interval $(-1,1)$.
The diffusion coefficient $k(x, \mu)$ in (\ref{mix-eq}) is defined as follows,
\[
k(x,\mu)=a(x,\mu)k_{hc},
\]
where $k_{hc}=\frac{10^4}{\kappa_1}$ and $\kappa_1$ is depicted in Figure \ref{fig-exam} (left), $a(x,\mu)$ take the first thirteen terms of $a(x,\mu)$, i.e., $N=12$ in KLE (\ref{KLE-TRUNC1}).
We note that $k^{-1}$ is not affine with respect to the parameter $\mu$. Here we apply STAOMP  to  $k^{-1}$ to get affine expression and achieve  offline-online computation decomposition. The reference solutions are computed by the mixed FEM on $60\times60$ fine grid. We compute the mixed GMsFE basis on  $10\times10$ coarse grid.

For the example, $500$ parameter samples are selected for the training set $\Xi_{\text{train}}$. Then we use the greedy algorithm (Algorithm \ref{algorithm-Greedy}) to select $|\Xi_{op}|=40$ optimal samples from the training set $\Xi_{\text{train}}$, and take $l_i(\mu)=5$ for each coarse block $\omega_i$ and each $\mu\in \Xi_{op}$, so the number of the snapshots is $n_{\text{snap}}^i=200$ for each coarse block $\omega_i$. Then we use RmGMsB-GBOCV to get the reduced multiscae basis  and apply LSMOS and STAOMP to the reduced model representation.
For LSMOS, we randomly choose $2000$ parameter samples for the snapshots, which are computed  by RmGMsB-GBOCV model with $7$ local basis functions.
 For LSMOS, we take the first $M$ KLE terms in (\ref{KLE Snapshots}). Here we take $M=3,4,5$, respectively. They are corresponding to the third, forth and fifth row in Table \ref{tab2-exam2}.
   Then use Legendre polynomials with total degree up to $N_g=4$ to approximate the KLE coefficients $\{\zeta_i(\mu)\}_{i=1}^{M}$, and the total number of the Legendre polynomial basis functions is $M_g=1820$.
   For STAOMP, we use POD to construct $6$ optimal global basis functions for spatial space, and choose
    $70$ parameter samples and $100$ spatial coordinates to get the reduced model representation by  Algorithm \ref{algorithm-STAOMP}, where the tolerance error $\varepsilon_\text{on}=1\times 10^{-4}$ for both velocity and pressure. To assess the approximation by the model reduction methods,
     we randomly choose $1000$ samples and compute the average relative errors defined as (\ref{errors-velocity}), (\ref{errors-pressure}). The results are listed in Table \ref{tab2-exam2} along with the average online CPU time per sample. In the table, $\varepsilon_v^1$ and $\varepsilon_p^1$ are the relative errors  based on the  reference solutions computed by mixed FEM on fine grid, while $\varepsilon_v^2$ and $\varepsilon_p^2$  are the relative errors based
      on the reference solution computed by  RmGMsB-GBOCV. From Table \ref{tab2-exam2}, we can see: (1) from the forth and sixth column of the table, we can see that the approximation  by the three approaches is almost identical, which means that  RmGMsB-GBOCV model can be accurately expressed  by the reduced model representation using LSMOS or STAOMP; (2) the average CPU time per sample by STAOMP  is the smallest among all of the approaches.   The numerical example shows that  the STAOMP  method achieves a good trade-off in both approximation accuracy and computation efficiency.
    We note that  the average online CPU time per sample is $21.4444$ seconds using  mixed FEM on fine grid. This time is much larger than the reduced model methods.

\begin{table}[hbtp]
\small
\centering
\caption{ Comparison of relative mean errors  and the average online CPU time for different approaches. $\text{Mt}_v$ is the number of terms of reduced model representation  for velocity, and $\text{Mt}_p$ is the number of terms of reduced model representation for pressure}
%\vspace*{2pt}
\begin{tabular}{c|c|c|c|c|c|c|c}
%\hline
\Xhline{1pt}
  Strategies & $\text{Mt}_v$ & $\text{Mt}_p$ & $\varepsilon_v^1$ &$\varepsilon_v^2$ &$\varepsilon_p^1$ &$\varepsilon_p^2$ & CPU time\\
  \hline
\multirow{1}{*}{RmGMsB}
               & & & $0.63 \times10^{-2}$ & &$2.76 \times10^{-2}$ &  & $0.3841 s$\\
 \hline
  \multirow{3}{*}{LSMOS}
               &  5460 &  5460 & $0.64 \times10^{-2}$ &$7.52\times10^{-5} $&$2.77\times10^{-2} $&$2.50\times10^{-3}$& $1.02\times10^{-2} s$\\
               \cline{2-8}
               &7280& 7280 & $0.63 \times10^{-2}$ &$6.27\times10^{-6}$&$2.76 \times10^{-2}$& $1.70\times10^{-3} $&$1.68\times10^{-2} s$\\
                \cline{2-8}
               & 9100& 9100 & $0.63 \times10^{-2}$ &$5.98 \times10^{-7}  $&$2.76 \times10^{-2} $&$6.60\times10^{-4} $& $1.84\times10^{-2} s$\\
 \hline
 \multirow{1}{*}{STAOMP}  &36 & 40 & $0.63 \times10^{-2}$ & $9.45\times10^{-6}$ &$2.76\times10^{-2}$& $1.40\times10^{-3} $ &  $2.02\times10^{-4} s$\\
\Xhline{1pt}
\end{tabular}
\label{tab2-exam2}
\end{table}

%%%%%%%%%%%%%%%%%%%%%%%%%%%%%%%%%%%%%%%%%%%%%%%%%%%%%%%%%%%%%%%%%%%%%%%%%%%

\subsection{Example of two-phase flow in random porous media}
\label{Numerical-result3}

In this numerical example, we consider the following two-phase flow problem  (in the
absence of gravity and capillary effects)  with zero Neumann boundary condition
\begin{eqnarray}
\label{two-phase-1}
\begin{cases}
\begin{split}
&v=\eta(S) k(x,\mu) \nabla p ~~\text{in} ~~D, \\
&-\nabla \cdot v=q ~~\text{in} ~~D,\\
&\frac{\partial S}{\partial t}+\nabla \cdot \big(f_w(S) v\big)=q_s  ~~\text{in}~~ D,
\end{split}
\end{cases}
\end{eqnarray}
where the total mobility $\eta(S)$ is defined  by $\eta(S)=S^2/\mu_w+(1-S)^2/\mu_o$ and $v$ refers to total velocity. Here $\mu_w/\mu_o=0.1$ is the ratio
between viscosity of water and oil, the fractional flow function  $f_w(S)$ is given by
\begin{eqnarray*}
\label{flux}
f_w(S)=\frac{S^2}{S^2+\mu_w/\mu_o(1-S)^2},
\end{eqnarray*}
 where $S$ is the saturation of water and $q$ is the source term for pressure equation and the source term for the saturation equation is expressed by
\[
q_s=\max(q,0)+f_w(S)\min(q,0).
\]
We consider a logarithmic random field,
\[
k(x,\mu):=\exp\big(\kappa_2+a(x,\mu)\big),
\]
where $\kappa_2$ is depicted in Figure \ref{fig-exam} (right), and $a(x,\mu)$ is a stochastic field, which is characterized by a two point exponential covariance function $\text{cov}[a]$ in (\ref{covariance function}).  Here the variance $\sigma^2=1$, correlation length $l_x=l_y=0.02$. The random field  $a(x, \mu)$ is obtained by  truncating  Karhunen-Lo\`{e}ve expansion   (\ref{KLE-TRUNC1}),
 where $E[a]=0$ and the random vector $ \mu:=(\mu_1, \mu_2, ...,\mu_N)\in \mathbb{R}^{N}$ and $\mu_i\sim N(0, 1)$  ($i=1,\cdots, N$), i.e, normal distribution with zero mean and unit variance.  We truncate the KLE (\ref{KLE-TRUNC1}) after the first $20$ terms to represent the random field $a(x, \mu)$.
  To  fulfil   offline-online computation decomposition,    we use EIM  for $k^{-1}(x, \mu)$ and  get an affine expression.

  The reference solution is computed on a $56 \times 56$ fine grid.  The mixed GMsFEM is performed on a $7 \times 7$  coarse grid.
   We solve the flow equation on the coarse grid using RmGMsB-GBOCV first and  then reconstruct the fine-scale velocity field as a superposition of the GMsFE basis functions. The reconstructed velocity field is used to solve the saturation equation with a finite volume method on the fine grid. We use an Implicit Pressure Explicit
Saturation (IMPES) formulation for the coupled system (\ref{two-phase-1}).

   Let $S_i^n$ be the saturation at the $i_{th}$ fine element $K_i^h$ at time $t_n=t_0+n \Delta t$, where $t_0$ is the initial time. Then the saturation equation can be discretized by
\[
|K_i^h|\frac{S_i^{n+1}-S_i^{n}}{\Delta t}+\int_{\partial K_i^h} f_w(\hat{S}^n)(v\cdot n)=|K_i^h|(q_s)_i,
\]
where $(q_s)_i$ is the average value of $q_w$ on $K_i^h$, and $\hat{S}^n$ is the upwind flux. For numerical test, we take $q$ to be zero except for the top-left and bottom-right coarse-grid elements, where $q$ takes the values 1 and $-1$, respectively. This is a traditional
two-spot problem,
in which the water is injected at the top-left corner and oil is
produced at the bottom-right corner.
Moreover, we set the initial staturation to be zero.

For the pressure equation in (\ref{two-phase-1}), we construct the multiscale finite element space $V_H^{M_C,CV}$ by RmGMsB-GBOCV at  initial time, and use it for all the time levels. We randomly choose $n_t=|\Xi_{\text{train}}|=200$ parameter samples and select $|\Xi_{op}|=40$ optimal samples by greedy algorithm to compute the snapshots. For each coarse block $\omega_i$ and each $\mu\in \Xi_{op}$, we take $l_i(\mu)=7$ local basis functions to construct snapshots ($n_{\text{snap}}^i=\sum_{j=1}^{n_t}l_i(\mu_j)=280$  for each coarse block $\omega_i$) for the reduced GMsFE basis.  Then we employ RmGMsB-GBOCV to simulate the flow equation and get the following weak formulation,
\begin{eqnarray}
\label{two-phase-dis}
\begin{cases}
\begin{split}
&a(v_H^N(\mu),u; \mu)-b(u, p_H^N(\mu); \mu)=0 ~~~\forall~ u \in V_H^{M_C,CV}, \\
&b(v_H^N(\mu), q; \mu)=l(q) ~~~\forall ~q \in Q_H,\\
&|K_i^h|\frac{S_i^{n+1}-S_i^{n}}{\Delta t}+\int_{\partial K_i^h} f_w(\hat{S}^n)(v\cdot n)=|K_i^h|(q_s)_i.
\end{split}
\end{cases}
\end{eqnarray}

\begin{figure}[htbp]
\centering
  \includegraphics[width=2.1in, height=1.9in]{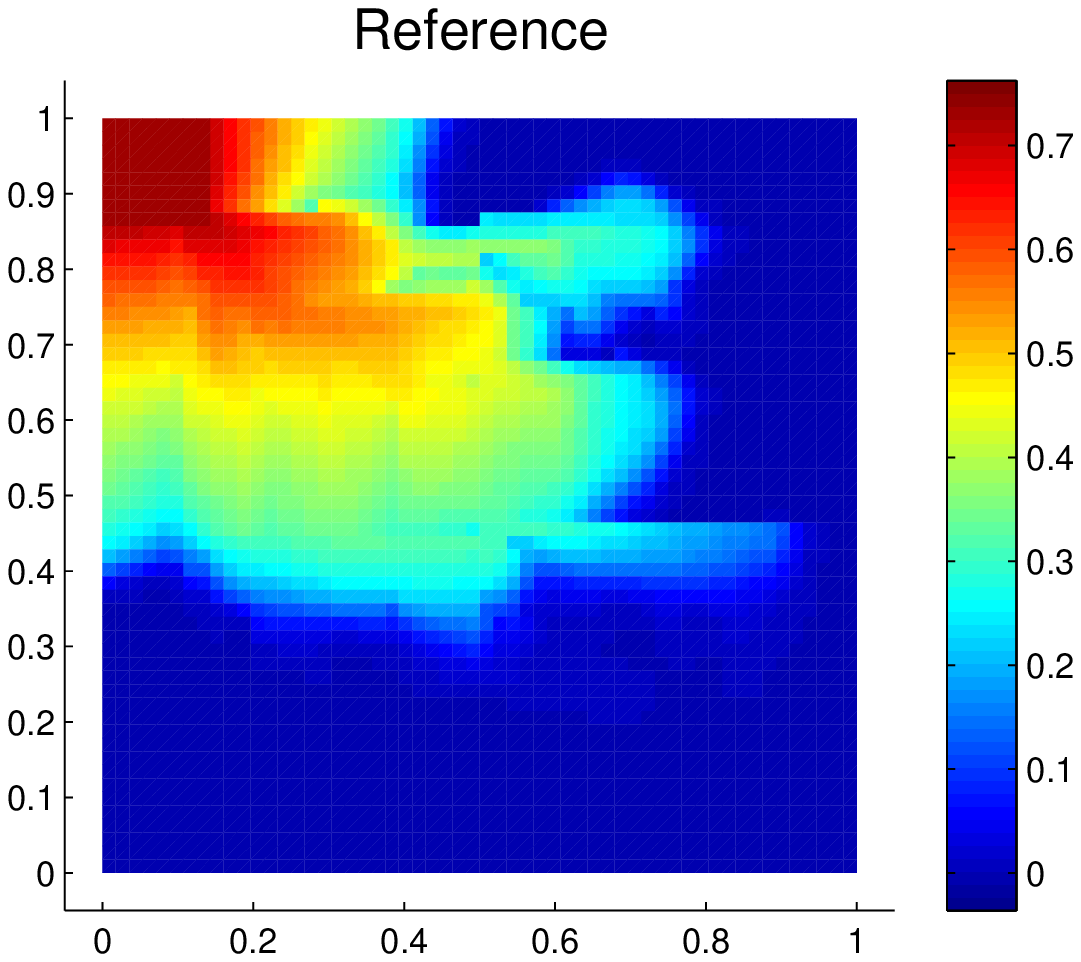}
  \includegraphics[width=2.1in, height=1.9in]{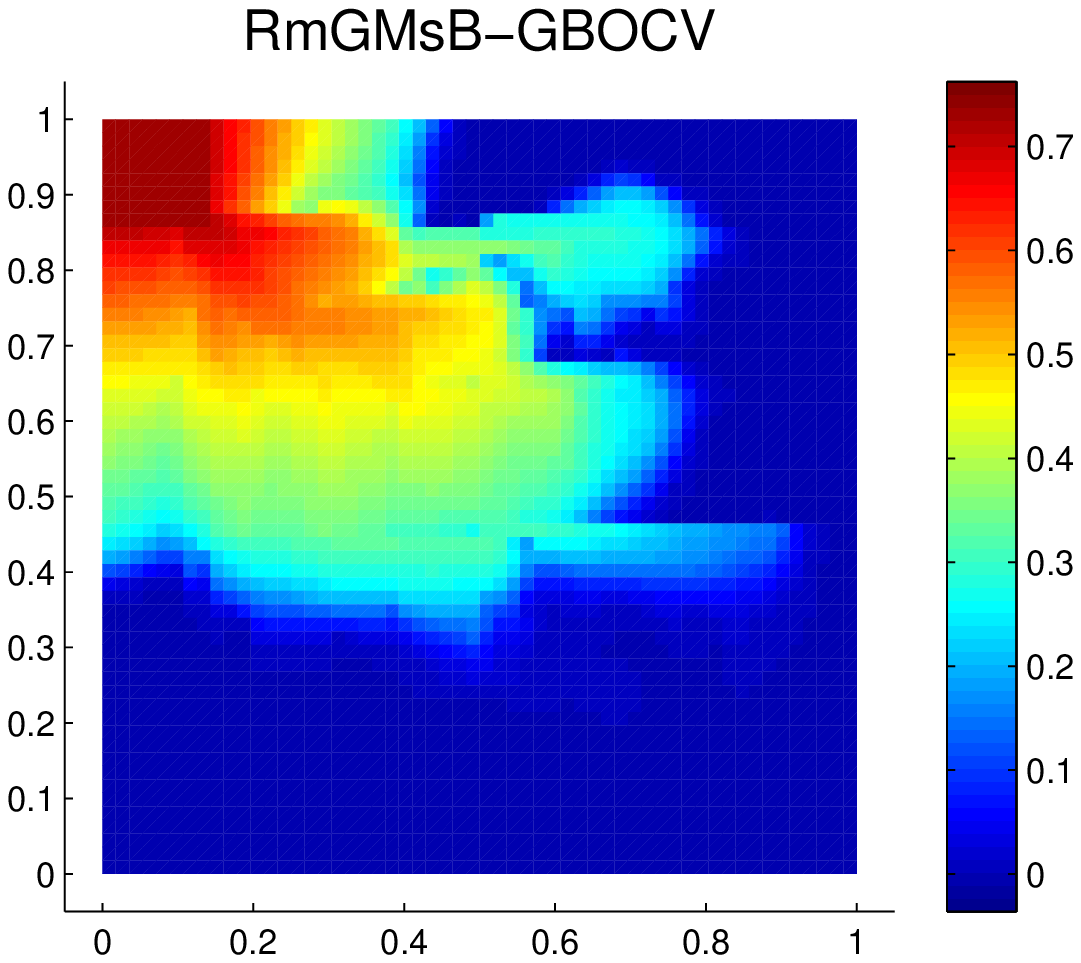}
  \includegraphics[width=2.1in, height=1.9in]{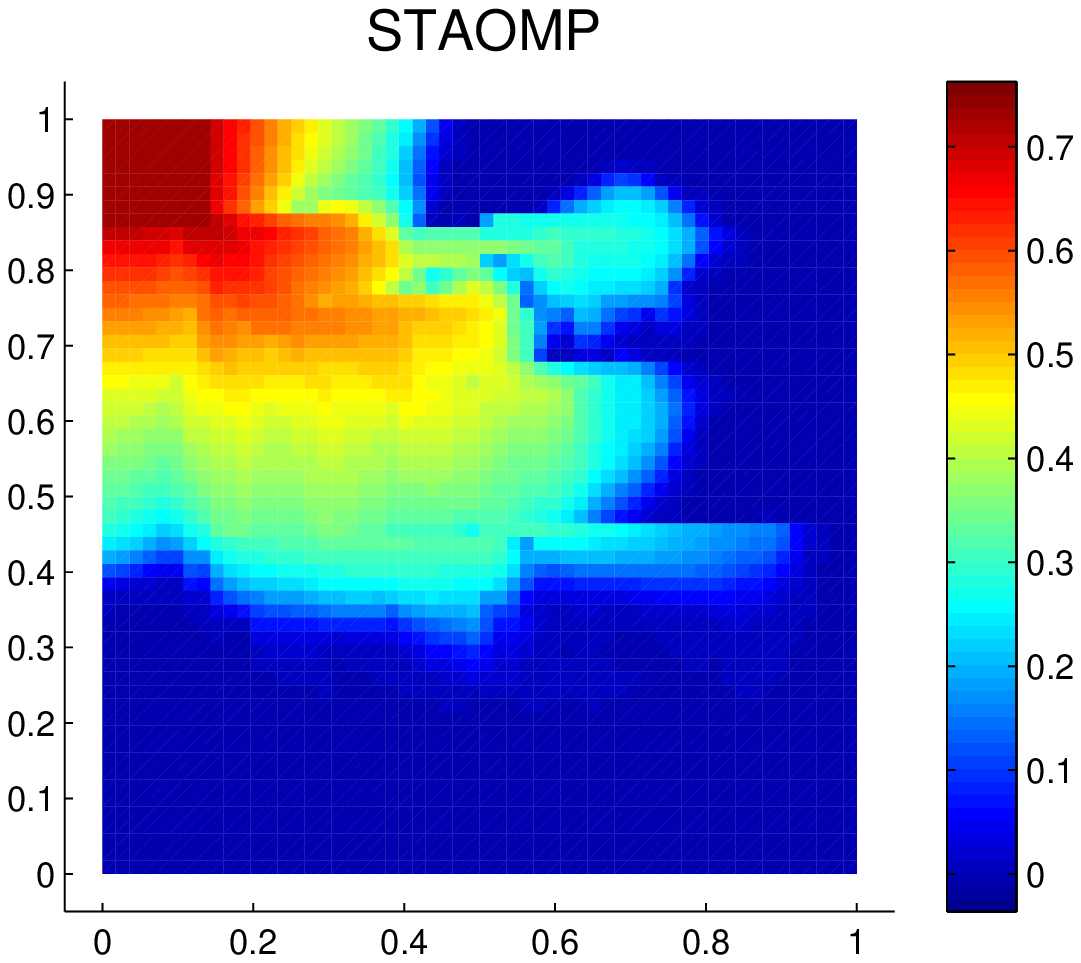}\\
  \includegraphics[width=2.1in, height=1.9in]{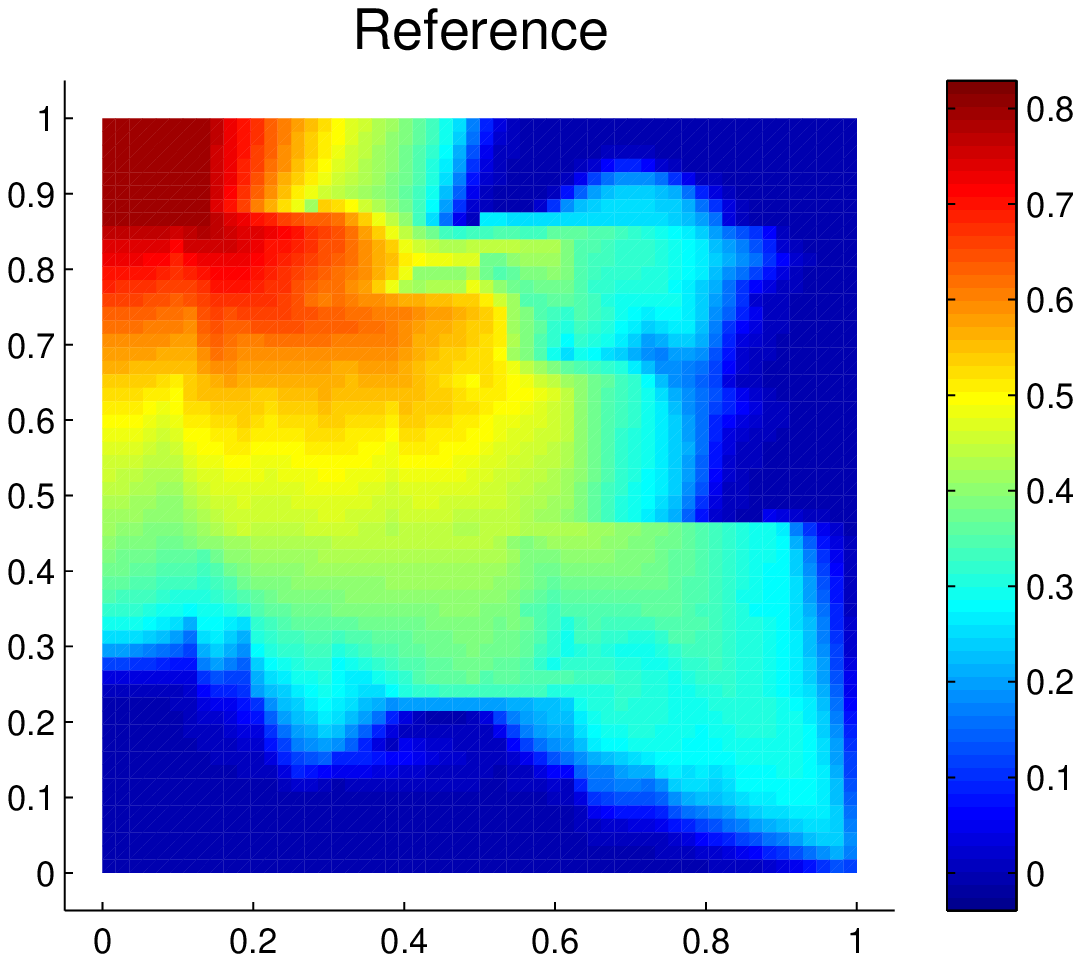}
  \includegraphics[width=2.1in, height=1.9in]{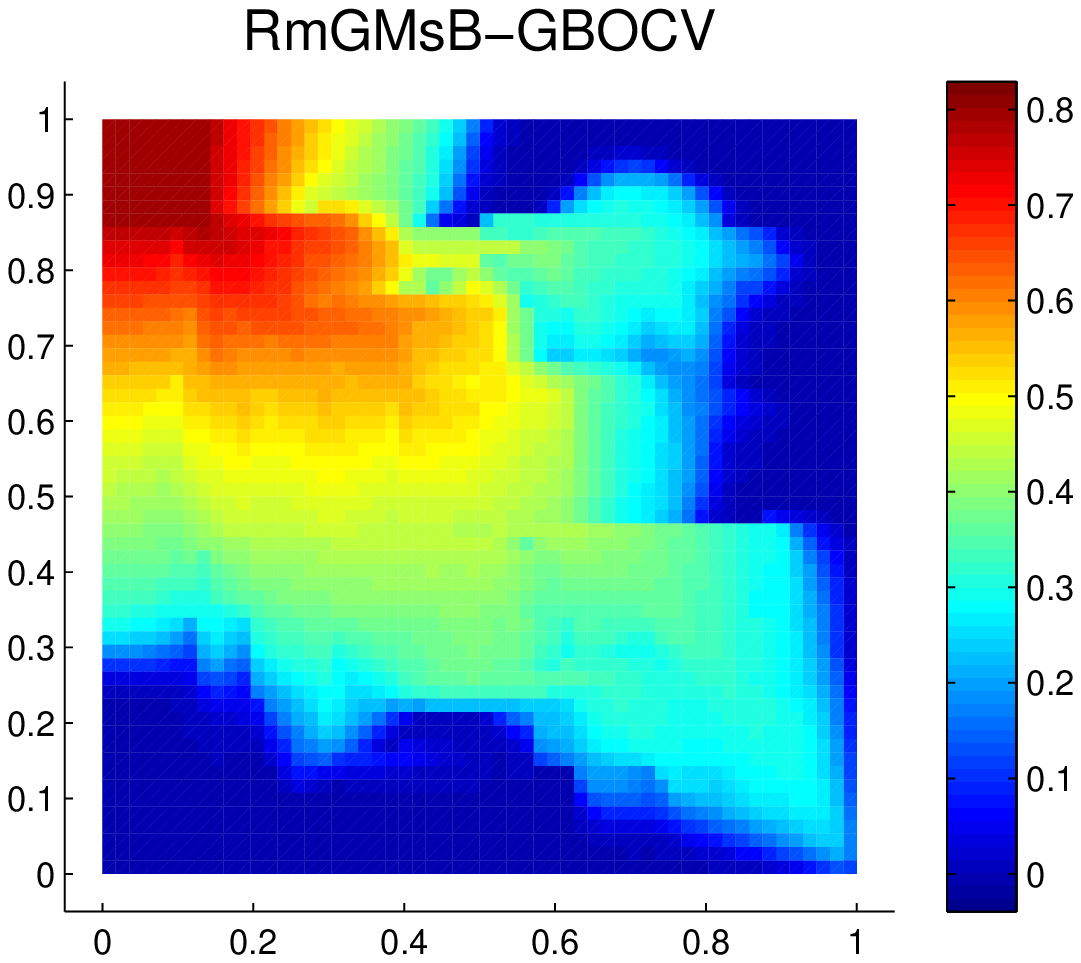}
  \includegraphics[width=2.1in, height=1.9in]{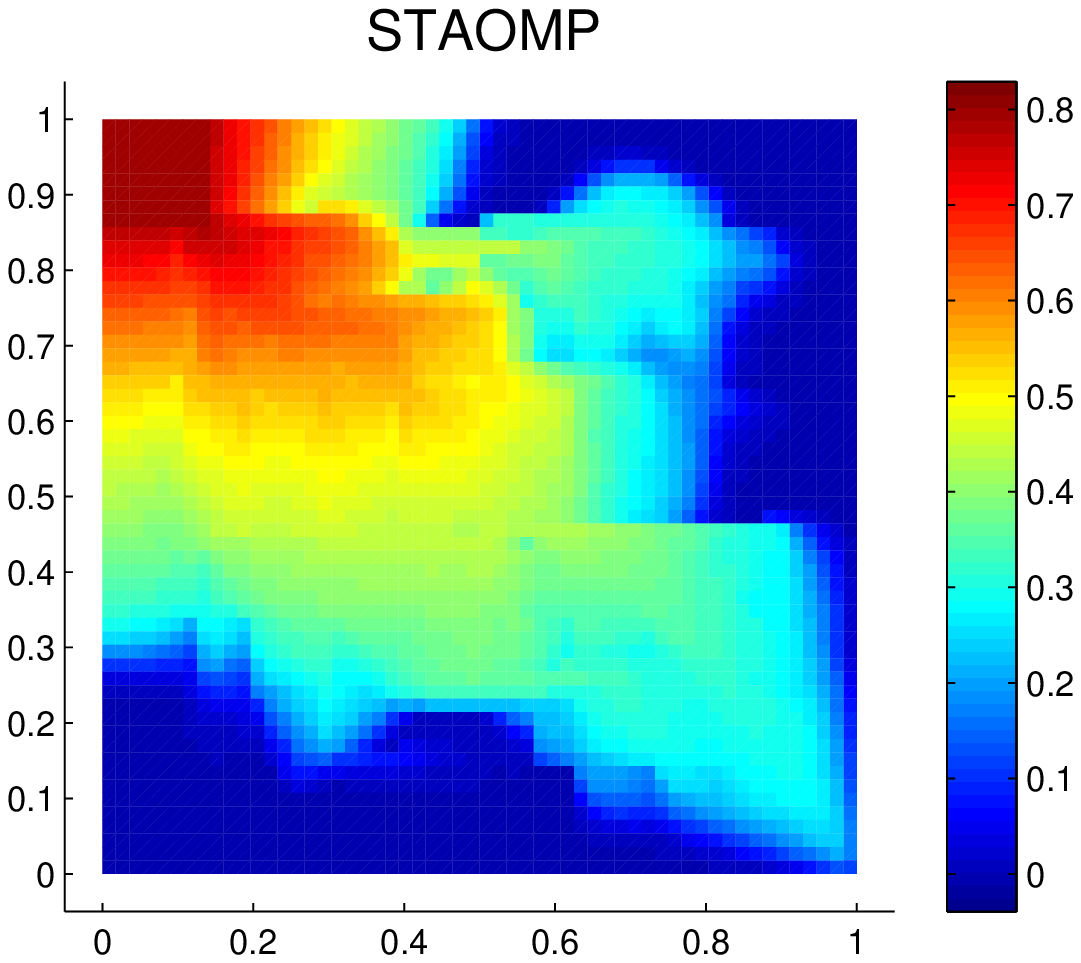}
  \caption{The mean of water saturation profiles at $t=600$ (the first row) and  $t=1000$ (the second row),  $56\times 56$ fine grid, and $7\times 7$ coarse grid, the number of local multiscale basis functions is $4$.}
  \label{fig1-exam3}
\end{figure}

\begin{figure}[htbp]
\centering
  \includegraphics[width=2.1in, height=1.9in]{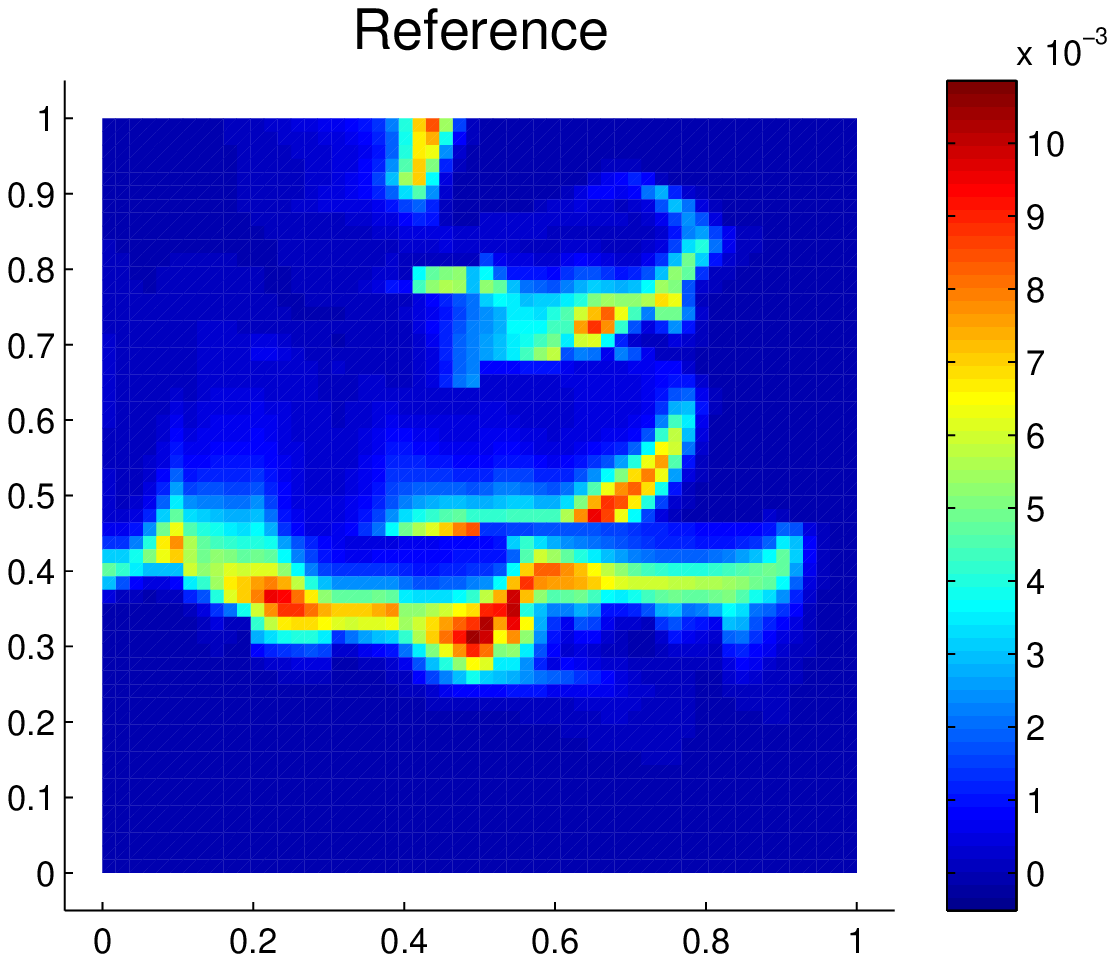}
  \includegraphics[width=2.1in, height=1.9in]{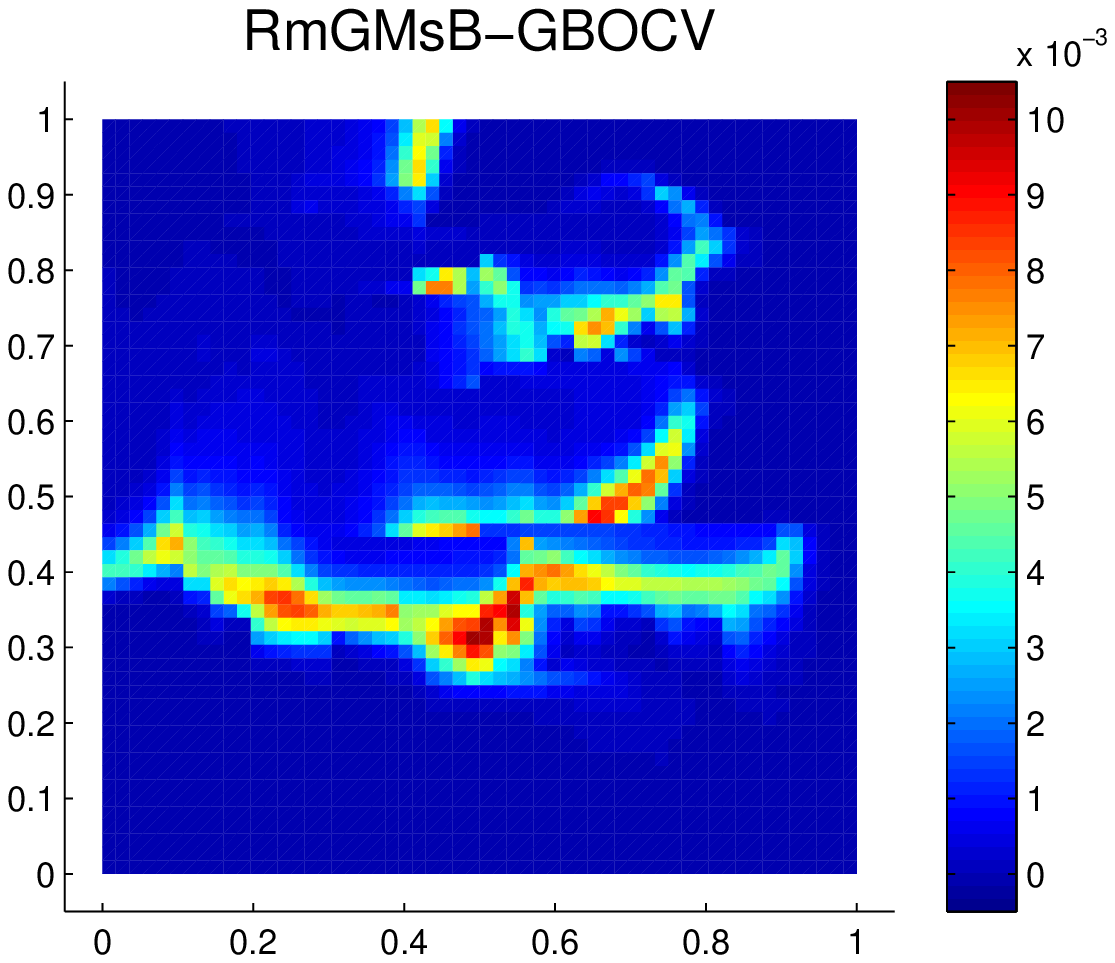}
  \includegraphics[width=2.1in, height=1.9in]{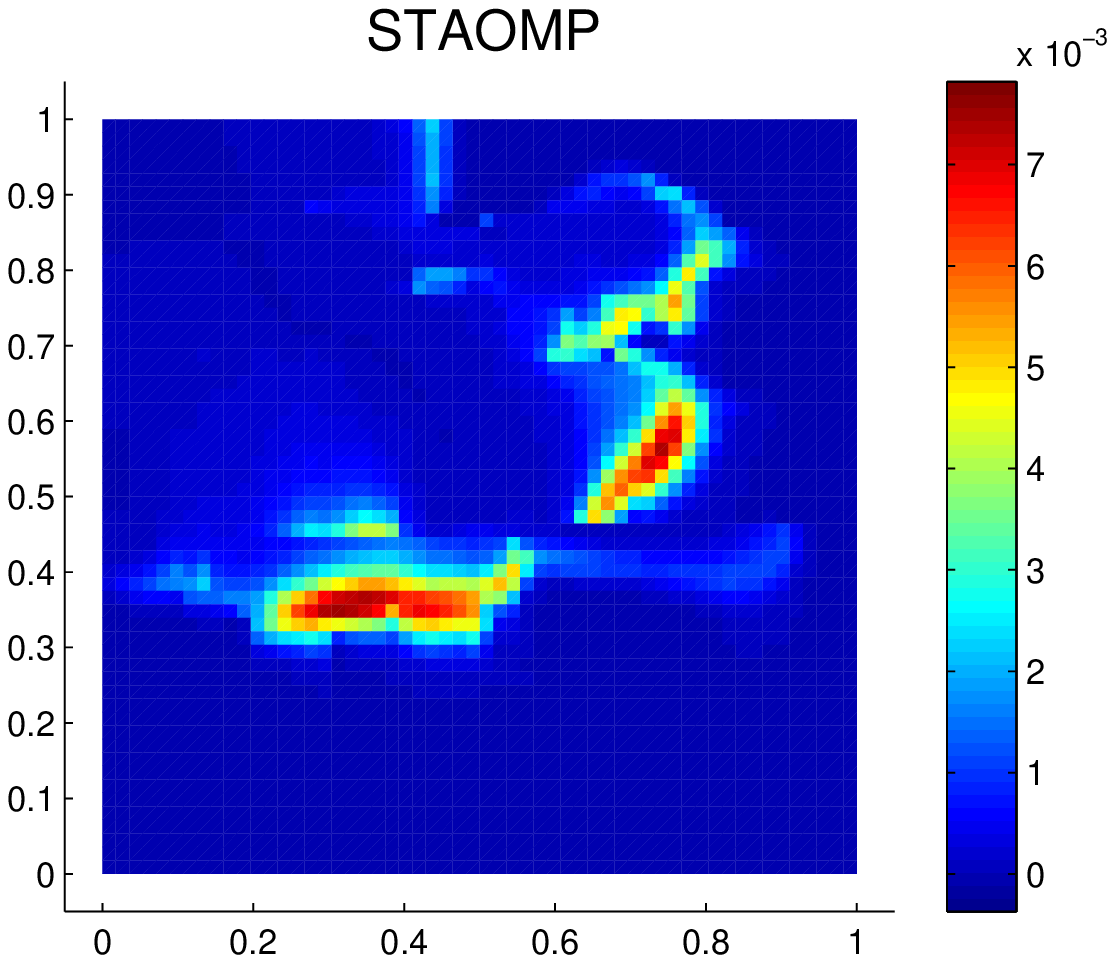}\\
  \includegraphics[width=2.1in, height=1.9in]{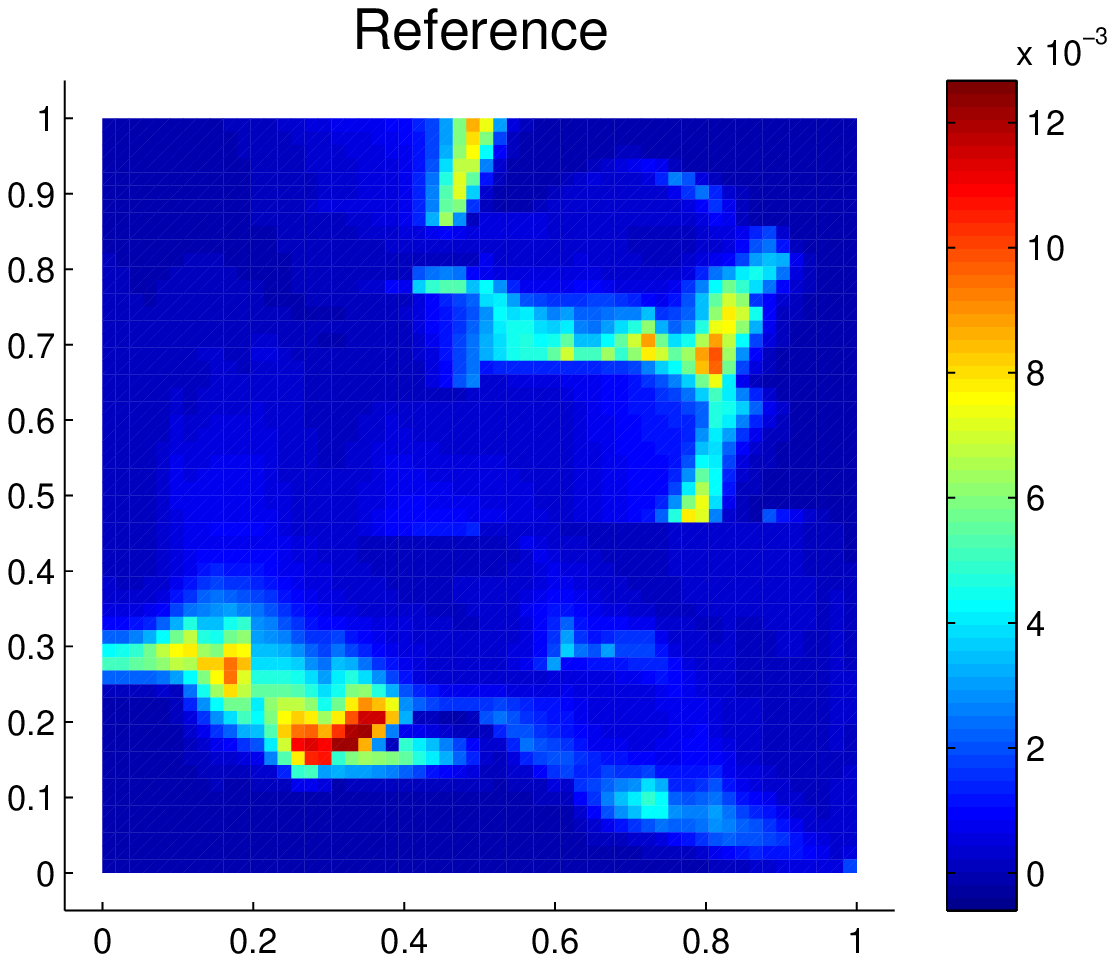}
  \includegraphics[width=2.1in, height=1.9in]{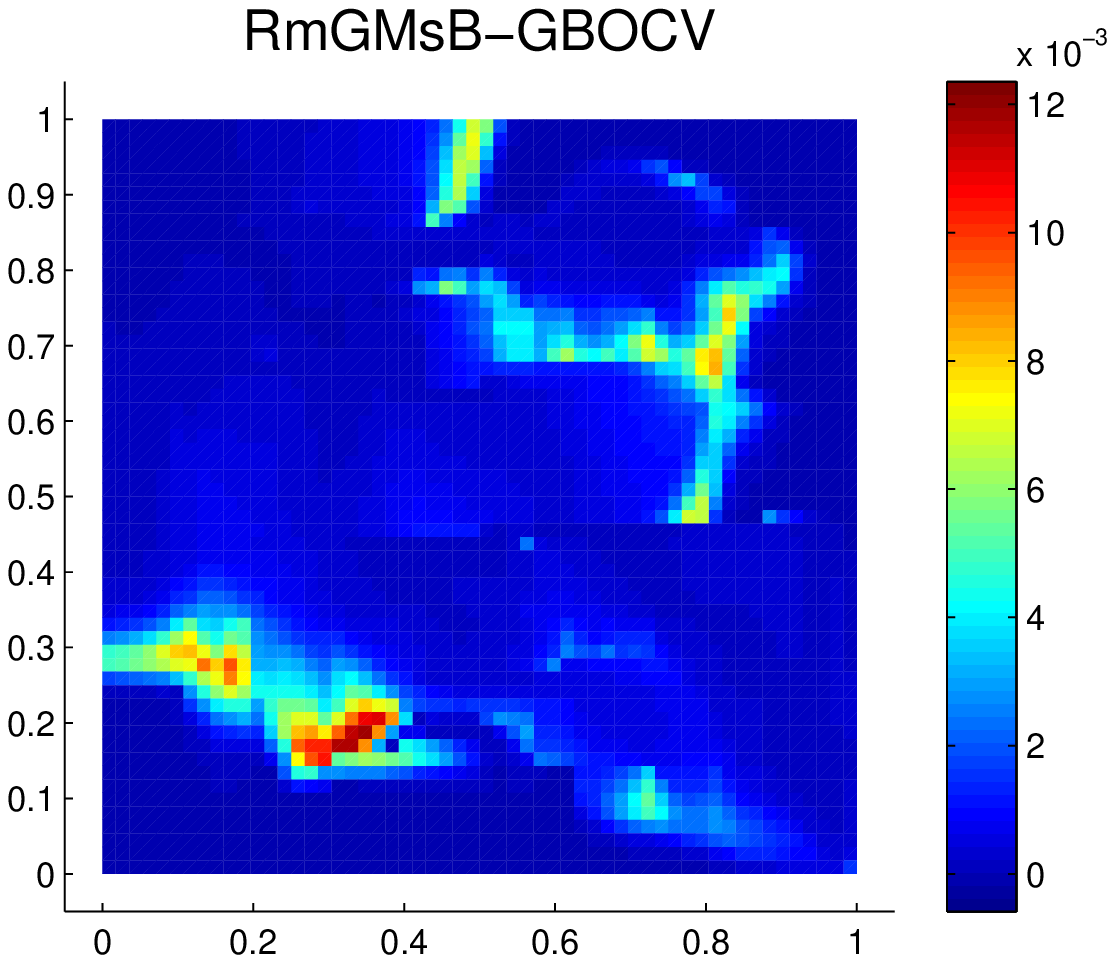}
  \includegraphics[width=2.1in, height=1.9in]{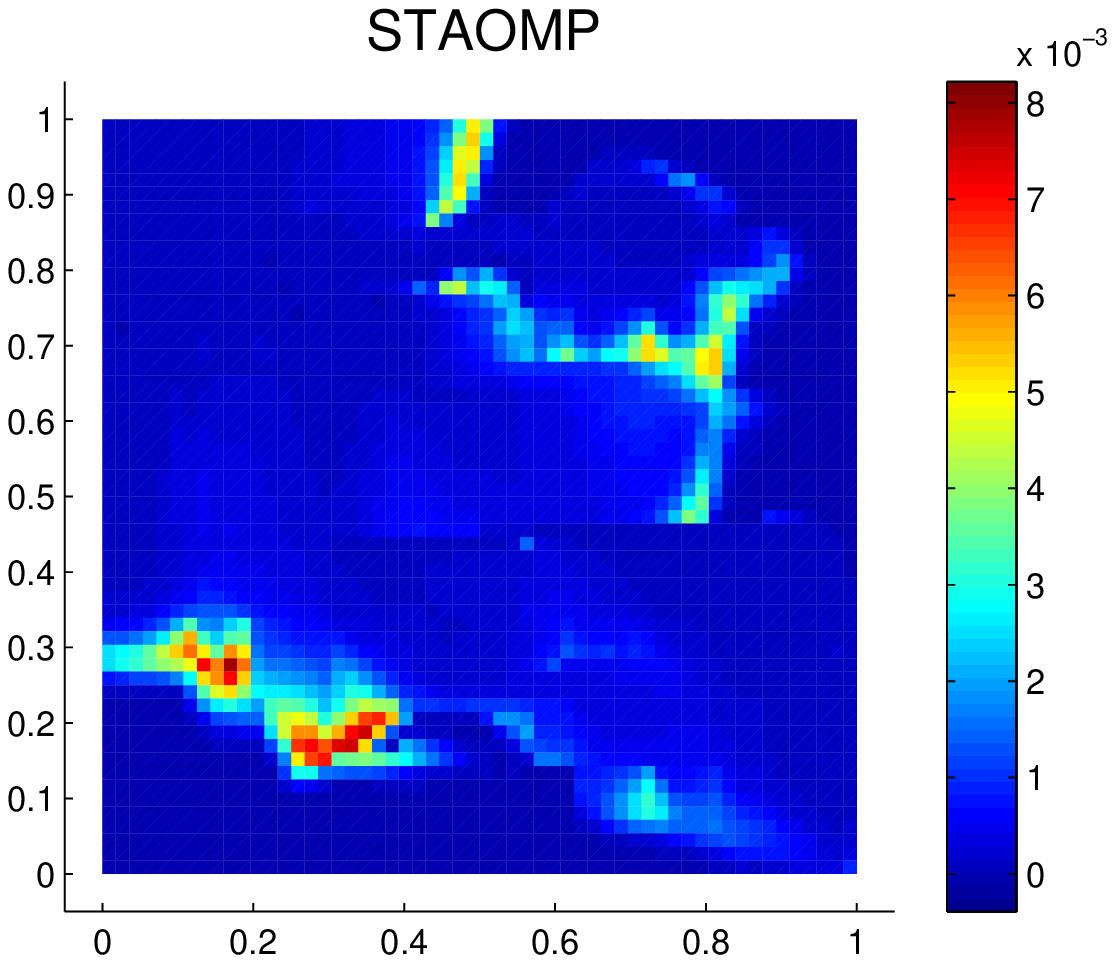}
  \caption{The variance of water saturation profiles at $t=600$ (the first row) and  $t=1000$ (the second row),  $56\times 56$ fine grid,  $7\times 7$ coarse grid, the number of local multiscale basis functions  is $4$.}
  \label{fig2-exam3}
\end{figure}

%\begin{figure}[htbp]
%\centering
%  \includegraphics[width=4.5 in, height=2.5in]{exp3/barcpu}
%  \caption{Comparison of the average online CPU time for water saturation per sample at $t=500$, $t=1000$, $t=2000$ and $t=1500$ for different approaches, the number of %local multiscale basis functions is $4$.}
%  \label{fig3b-exam3}
%\end{figure}

\begin{table}[hbtp]
\small
\centering
\caption{ Comparison of the average online CPU time for water saturation per sample at different time levels for the two approaches: STAOMP and RmGMsB-GBOCV, the number of local multiscale basis functions is $4$.}
%\vspace*{2pt}
\begin{tabular}{c|c|c|c|c}
%\hline
\Xhline{1pt}
  Strategies & $t=500$ & $t=1000$ & $t=1500$ & $t=2000$\\
  \hline
   STAOMP  & $1.60\times 10^{-3} s$ & $1.50\times10^{-3} s$ & $1.40\times10^{-3} s$ & $1.50\times10^{-3} s$\\
 \hline
 RmGMsB-GBOCV & $1.7837$ & $5.0227$ & $6.6340$ & $7.8311$\\
\Xhline{1pt}
\end{tabular}
\label{table3b-exam3}
\end{table}

Water saturation is an important quantity for the model.  We utilize  STAOMP  to represent  the water saturation $S(x,t,\mu)$ of the reduced model (\ref{two-phase-dis}).
 To this end,  we choose the  set of parameter samples  $\Xi_{t}=\{\mu_i\}_{i=1}^{600}$ and apply POD to the snapshots associated with $\Xi_{t}$ and get the optimal  global basis  functions $\{\hat{S_i}(x,t)\}_{i=1}^{4}$.
To approximate random parameter space, we use Hermite polynomial basis functions with total degree up to $N_g=3$ and denote the basis set by
  $\{p_i(\mu)\}_{i=1}^{1771}$.
   A finite dimensional approximation space  for $\{S(x,t,\mu_i); \mu_i\in \Gamma\}$ is then obtained by
\[
\textbf{S}_{4\times 1771}:=\text{span}\{p_i(\mu)\hat{S_j}(x,t): 1\leq i\leq 1771,1\leq j \leq 4\}.
\]
%To simplify the notation, we use the following single-index notation
%\[
%\textbf{S}_{7084}=\{W(x,t,\mu)=\sum_{i=1}^{7084}w_i\Psi_i(x,t,\mu); ~w_i\in \mathbb{R}\}.
%\]
To get a model reduction representation using  STAOPM, we randomly choose $600$ parameter samples, $64$ spatial points and 18 time levels  for snapshots, and then use Algorithm \ref{algorithm-STAOMP} to construct
 a sparse representation $S(x,t,\mu)\approx \sum_{i=1}^{\text{Mt}} \textbf{c}(\emph{I}(i)) \Psi_{\emph{I}(i)}(x,t,\mu)$, where $\text{Mt}=160$,
 $I=\{1,\cdots,1771\}\times \{1,\cdots,4\}$ and $\Psi_i(t,\mu)=p_{i_1}(\mu)\hat{S_{i_2}}(x,t)$.

 To calculate the approximation errors by RmGMsB-GBOCV and  STAOPM,  we randomly choose $1000$ samples and compute the RmGMsB-GBOCV model (\ref{two-phase-dis}) and evaluate the reduced model representation by STAOMP  for each sample. The relative mean  errors for saturation $S$ is defined by
\begin{eqnarray}
\label{errors-saturation}
\varepsilon_s(t)=\frac{1}{N}\sum_{i=1}^N\frac{\|S(x,t,\mu_i)-S_{H}(x,t,\mu_i)\|_{L^1(D)}}{\|S(x,t,\mu_i)\|_{L^1(D)}},
\end{eqnarray}
where $N=1000$,  $S(x,t,\mu_i)$ is the solution of fine scale model, and   $S_{H}(x,t,\mu_i)$ is the solution of (\ref{two-phase-dis}) by RmGMsB-GBOCV or STAOPM.

\begin{figure}[htbp]
\centering
  \includegraphics[width=4.5 in, height=2.5in]{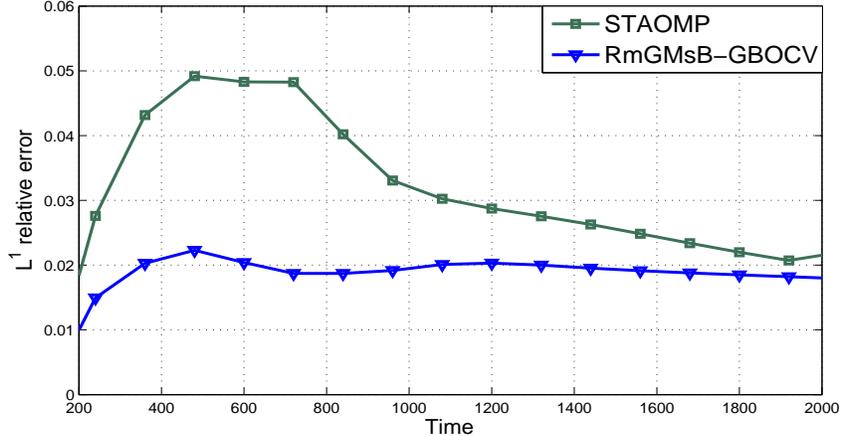}
  \caption{The average relative error for water saturation versus different time levels. the number of local multiscale basis functions is $4$.}
  \label{fig3-exam3}
\end{figure}

Figure \ref{fig1-exam3} shows the  mean of water saturations at $t=600$ and $t=1000$ for the different methods.  Here the reference solution is computed  by mixed FEM on fine grid. By the figure, we see that
the model reduction methods provide accurate approximation to the reference saturation in the average sense.
Figure \ref{fig2-exam3} demonstrates  the variance of water saturation at $t=600$ and $t=1000$  for the different methods. By the figure, we can find: (1) the largest variance occurs along the advancing water front;  (2)
there is no clear difference for the variance profiles  between the reference solution and  RmGMsB-GBOCV solution; (3) compared with reference solution, the STAOMP approach renders some clear difference for the variance approximation around the water front.
Table \ref{table3b-exam3} lists the average online CPU time per sample for the water saturation by the different approaches. We find that the STAOMP  method significantly improves computation efficiency especially for the later time levels. Moreover, the online CPU time of STAOMP method is almost independent of the time level.
This is an important advantage of the STAOMP  method for two-phase  flow simulation. Figure \ref{fig3-exam3} shows the average relative error for water saturation defined by (\ref{errors-saturation}) versus different time levels for the different methods. By the figure, we find that
%:(1) for  RmGMsB-GBOCV,  the relative error for water saturation become smaller as the number of local basis functions increases; (2)
the error by  STAOMP decays after the water break-through time and then tends to be stable.

\begin{figure}[htbp]
\centering
  \includegraphics[width=3in, height=2.5in]{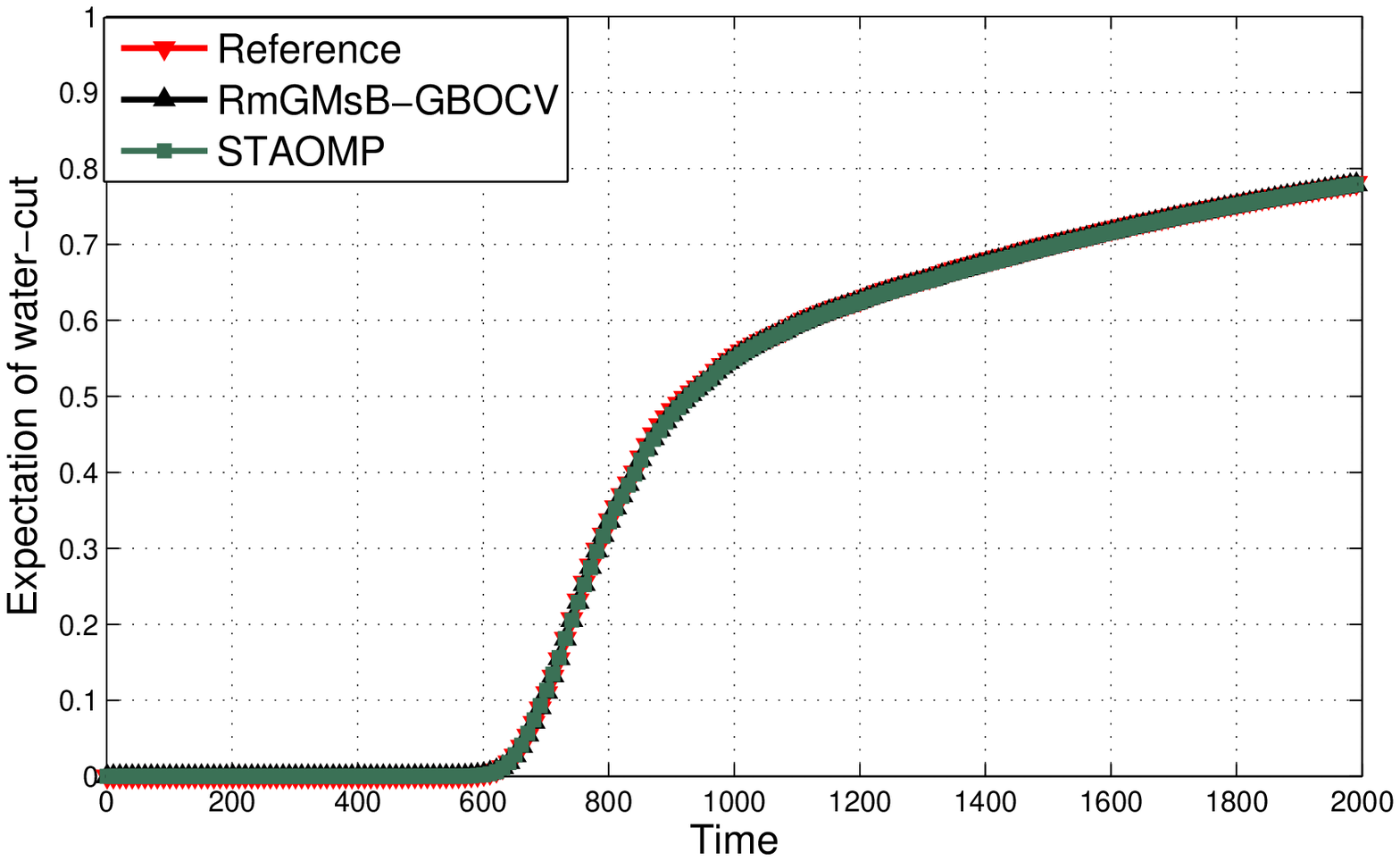}
  \includegraphics[width=3in, height=2.5in]{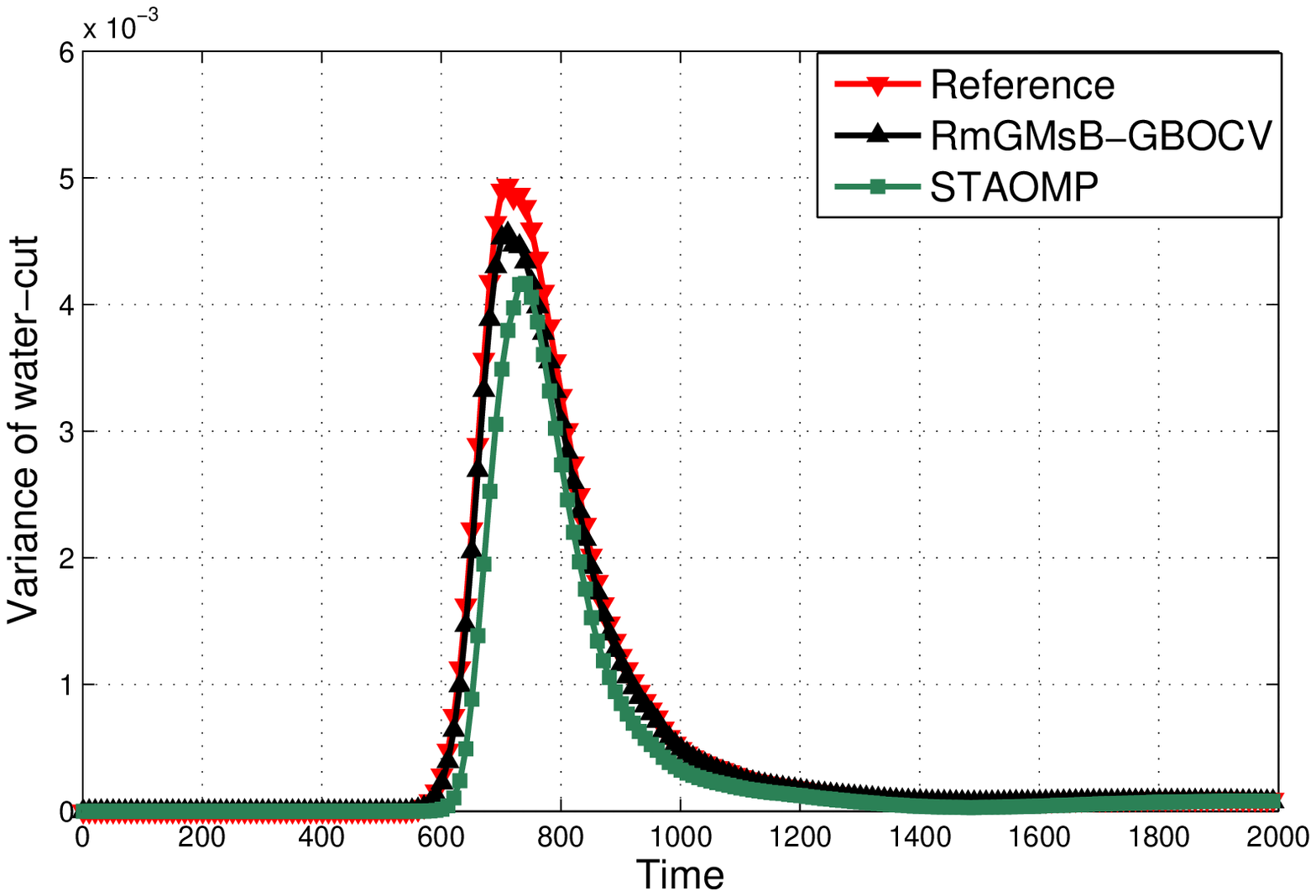}
  \caption{The expectation (left) and variance (right) of water-cut curves}
  \label{fig4-exam3}
\end{figure}

Water-cut is also an important quantity for the water-oil two-phase flow. The water-cut  is defined as the fraction of water in
the produced fluid and is given by $q_w/q_t$, where $q_t=q_o+q_w$, with $q_o$ and $q_w$ being the flow rates of oil
and water at the production edge of the model. In particular, $q_w=\int_{\partial D^{out}} f_w(S){v}\cdot { n} ds$,
$q_t=\int_{\partial \Omega^{out}} {v}\cdot { n} ds$,
where ${\partial D^{out}}$ is the out-flow boundary.

For the water-cut, we utilize  STAOMP  to approximate the water-cut of the reduced model (\ref{two-phase-dis}).
We use the set of parameter samples  $\Xi_{t}=\{\mu_i\}_{i=1}^{600}$ and apply POD to the snapshots associated with $\Xi_{t}$ to obtain $6$ optimal global basis functions $\{\hat{w_i}(t)\}_{i=1}^{6}$. We use the same polynomial basis functions as the water saturation, and get a finite dimensional approximation space
\[
\textbf{W}_{6\times 1771}:=\text{span}\big\{p_i(\mu)\hat{w_j}(t): 1\leq i\leq 1771,1\leq j \leq 6\big\}.
\]
To get a model reduction representation using  STAOPM, we randomly choose $600$ parameter samples and $10$ time levels for snapshots, and then use Algorithm \ref{algorithm-STAOMP} to construct
a sparse representation $W(t,\mu)\approx \sum_{i=1}^{\text{Mt}} \textbf{c}(\emph{I}(i)) \Psi_{\emph{I}(i)}(t,\mu)$ with $\text{Mt}=195$,
where $I=\{1,\cdots,1771\}\times \{1,\cdots,6\}$ and $\Psi_i(t,\mu)=p_{i_1}(\mu)\hat{w_{i_2}}(t)$.

Figure \ref{fig4-exam3} shows the expectation (left) and variance (right) of the water-cut curves with the different approaches. The figure shows that the expectation curves of the water-cut are all nearly identical. This implies a good mean approximation for water-cut using  RmGMsB-GBOCV and STAOMP; The variance of water-cut rises rapidly at the break-through time, and then decays fast after a while. The error of water-cut variance by
 STAOMP is slightly larger around the water break-through  time than other time instances.

\section{Conclusions}
\label{ssec:Conclusions}
In the paper, we presented  two variable-separation representations for random fields: LSMOS and STAOMP.  The LSMOS is devoted to constructing a KL expansion and building a relation  between the random inputs  and the stochastic basis functions based on  least-squares methods.  The STAOMP provides  a tensor product for the random field and achieves a sparse approximation by using orthogonal-matching-pursuit method.
  The two methods can be applied to get a variable-separation approximation  for any generic functions with random inputs.
For the multiscale  problems with high-dimensional random inputs, the computation would be prohibitive  if we  directly solve the problems for a many-query situation. To overcome the difficulty,
we developed reduced mixed GMsFE basis methods to improve the computation efficiency. To this end, we have chosen  a few optimal parameters from the parameter space by a greedy algorithm. Then we use the optimization methods  BOCV and POD  to obtain reduced multiscale basis functions and get a reduced order model. In  the reduced multiscale basis methods, the whole computation admits offline-online decomposition. In the online phase,  a reduced model
is solved for each parameter sample.
Although the offline computation may be expensive, the online computation is  efficient. We have carefully compared the performance of the different reduced mixed GMsFE basis methods, and found that BOCV  may lead to the reduced model with better accuracy   and robustness  than POD does.  The dimension of the space spanned by the set of reduced  multiscale basis functions is much smaller than the dimension of the original full order model, but it depends on the size of coarse grid.  When the size of coarse grid is large, the computation of the reduced order model may be still expensive.
To further improve the efficiency of the online computation,   we construct a sparse representation for model's outputs by combining the reduced mixed GMsFE basis methods and  STAOMP.
This is very desirable for predicting the model's outputs for arbitrary  parameter values and quantifying the uncertainty propagation. We applied the proposed reduced mixed GMsFE basis methods and the model's sparse representation method to a few numerical models with multiscale and random inputs. Careful numerical analysis is carried out for these computational models.

Although the proposed methods significantly reduces the online computation, it takes much effort in offline computation stage. In future, we will explore novel methods
to substantially decrease the offline computation burden as well.

\section*{Acknowledgments}
We would like to thank the helpful discussion with  Prof. Yalchin Efendiev for the work.    We acknowledge the support of Chinese NSF 11471107.

\end{document}